\documentclass[times,sort,3p]{elsarticle}

\usepackage[labelfont=bf]{caption}

\usepackage{amsmath,amsfonts,amssymb,amsthm,booktabs,color,epsfig,graphicx,hyperref,url}
\usepackage{xcolor}
\usepackage{enumitem}
\theoremstyle{plain}
\newtheorem{theorem}{Theorem}

\newtheorem{proposition}{Proposition}
\newtheorem{lemma}{Lemma}

\theoremstyle{definition}

\newtheorem{example}{Example}

\newenvironment{myproof}[1] {\paragraph{\textnormal{\textbf{Proof of {#1}}}}}{\hfill$\square$}

\usepackage[numbers]{natbib}
\def\mR{\mathbb{R}}
\def\tr{\mbox{tr}}
\def\cov{\mbox{cov}}
\def\var{\mbox{var}}
\def\sign{\mbox{sign}}

\def\diag{\mbox{diag}}

\def\defby{\stackrel{\mbox{\textrm{\tiny def}}}{=}}
\newcommand{\bSig}{\mbox{\boldmath $\Sigma$}}
\newcommand{\sbSig}{\boldmath{\scriptstyle \Sigma}}

\newcommand{\bI}{\mathbf I} 
\newcommand{\one}{\mathbf 1} 
\newcommand{\bzero}{\mathbf 0} 

\newcommand{\bH}{\mathbf{H}}

\newcommand{\E}{\mathbb E}

\newcommand{\bS}{\mathbf S}
\newcommand{\X}{\mathbf X}

\newcommand{\Z}{\mathbf Z}

\newcommand{\R}{\mathbf R}

\newcommand{\A}{\mathbf A}

\newcommand{\W}{\mathbf W}

\newcommand{\trans}{^\top}
\newcommand{\btau}{\boldsymbol \tau }
\newcommand{\brho}{\boldsymbol \rho }

\newcommand{\eqlsd}{\overset{\textrm{\tiny LSD}}{=\joinrel=}}

\begin{document}
	\begin{frontmatter}
	\title{Limiting spectral distribution of large dimensional Spearman's rank correlation matrices\tnoteref{mytitlenote}}
	
	\author{Zeyu Wu}
	\author{Cheng Wang\corref{mycorrespondingauthor}}
	\cortext[mycorrespondingauthor]{Corresponding author. Email address:  \url{chengwang@sjtu.edu.cn}}
	
	\address{School of Mathematical Sciences, MOE-LSC,\\ Shanghai Jiao Tong University,  Shanghai 200240, China.}


	\begin{abstract}
		In this paper, we study the empirical spectral distribution of Spearman's rank correlation matrices, under the assumption that the observations are independent and identically distributed random vectors and the features are correlated.    We show that the limiting spectral distribution is the generalized Mar\u{c}enko-Pastur law with the covariance matrix of the observation after standardized transformation.  With these results, we compare several classical covariance/correlation matrices including the sample covariance matrix, Pearson's correlation matrix, Kendall's correlation matrix and Spearman's correlation matrix.
	\end{abstract}
	
	
\end{frontmatter}

\section{Introduction}
Statistical inference for covariance or correlation matrices are fundamental  problems in high dimensional data analysis \cite{fan2016overview,cai2017global}. Recently, rank-based correlation matrices, e.g., Kendall's tau and Spearman's rho, have drawn increasing attention in a variety of fields.  Due to the robustness of non-parametric statistics, these rank-based correlations have appealing properties in theory and have the potential to solve problems for high dimensional data with complex structure such as heavy-tailed distributions \cite{hettmansperger2010robust}.

In high dimensional data analysis, \citet{liu2012high} and \citet{xue2012regularized}  firstly used rank-based correlation matrices to conduct sparse estimation for the covariance matrix and the precision matrix (i.e, the inverse of the covariance matrix). For Gaussian distribution, there exists an explicit relationship between rank-based correlation and Pearson's correlation. Exploiting this neat property, \citet{liu2009nonparanormal} extended the Gaussian graphical model to the non-paranormal distribution. Further, there are many works on hypothesis testing of covariance/correlation matrices, which are based on Kendall's tau correlation matrix $\btau_n$ or Spearman's rho correlation matrix $\brho_n$. In details,  \citet{bao2015Spearman} derived the asymptotic normal distributions of $\tr(\brho_n^k),~k\in \{2,3,\ldots\}$ and conducted several statistics for testing the complete independence. \citet{han2017distribution} studied the maximum norm of off-diagonal elements which converges weakly to Gumbel distribution and proposed test statistics based on  Kendall's tau and Spearman's rho. \citet{leung2018testing}  considered the test statistics based on the squared Frobenius norm and see also \citet{mao2017robust}.  


Technically,   \citet{liu2012high} and \citet{xue2012regularized} focused on the sparse estimation of Kendall's tau and Spearman's rho.  To guarantee the consistency,  they need to bound the maximum norm between the empirical and population rank-based correlation matrices, i.e., $\|\brho_n-\brho\|_\infty$ and $\|\btau_n-\btau\|_\infty$.  To control the maximum norm, the crucial issue is to bound the tail probability of each element uniformly. The main tools are the Hoeffding decomposition for U-statistics and some classical concentration inequalities such as Hoeffding and McDiarmid inequalities. As far as the hypothesis testing,  to derive the asymptotic distribution of the test statistics, we need more refined results. Under the null hypothesis of independence, all the features are independent and it is doable to derive the distribution of the test statistics. Under the alternative case, it is challenging to analyze Kendall's tau or Spearman's rho since rank variables depend on all the samples. Thus, the existing works \cite[e.g.,][]{bao2015Spearman,han2017distribution,mao2017robust,leung2018testing} only derived the distribution under the null case.

In random matrix theory,  the same challenges also exist. \citet{bai2008large} and \citet{bandeira2017marvcenko} studied the limiting spectral distribution of Spearman's rho  and Kendall's tau, respectively.  \citet{bao2015Spearman} and \citet{li2021central} considered the asymptotic distribution of the test statistics involved Spearman's rho  and Kendall's tau, respectively. \citet{bao2019tracy_sp} and \citet{bao2019tracy} derived the Tracy--Widom limits for Spearman's rho  and Kendall's tau, respectively.  All these works are based on the independent assumption. For the data with a common covariance structure, the data matrix after ranking does not have independent columns or rows anymore and in random matrix theory, it is challenging to analyze such matrices and see \citet{bai2008large} for more details. Recently, we derived the limiting spectral distribution (LSD) of Kendall's rank correlation matrix with dependence in \citet{li2021eigenvalues} and it is the first result on rank correlation matrices with dependence. In this paper, as a companion work, we study the limiting spectral distribution of Spearman's rank correlation matrices.  Under mild conditions, we show that the LSD follows the  generalized Mar\u{c}enko-Pastur law with a conditional covariance matrix. With these results, we can look insider into Spearman's correlation matrix, and also its connections with other covariance/correlation matrices. 

The rest of the paper is organized as follows. In Section \ref{sec2}, we revisit the Spearman's rank correlation matrix and provide some trivial results when the data is correlated. In Section \ref{sec3}, we formulate Spearman's correlation matrix in the form of U-statistics and we present the main results on LSD in Section \ref{sec4}. Finally, we provide discussions on several covariance/correlation matrices and also compare their LSDs through numerical simulations in Section \ref{sec5}. All technical proofs are relegated to Appendix.

\section{Spearman's rank correlation matrix} \label{sec2}
We first introduce some necessary notation. Throughout the paper,  $\bI_m$ is a $m \times m$ identity matrix and $\one_m$ is a $m \times 1$ vector with all ones. $\|\cdot\|_2$ denotes the Frobenius norm of a vector or matrix. $\|\cdot\|$ is the spectral norm and $\|\cdot\|_\infty$ is the element-wise maximum norm of a vector or matrix, respectively.

Given the  independent and identically distributed (i.i.d.)  observations $\X_1,\ldots, \X_n \in \mR^p$, we have the data matrix
\begin{align*}
	\X=
	\begin{pmatrix}
		\X_1\trans\\
		\vdots\\
		\X_n\trans
	\end{pmatrix}=\begin{pmatrix}
		x_{11} &\cdots & x_{1p}\\
		\vdots &\ddots& \vdots\\
		x_{n1}&\cdots & x_{np}
	\end{pmatrix}=
	(x_{ij})_{n \times p},
\end{align*}
where each row is an observation and each column is a feature. Transforming each column of the data matrix into the order statistics $(r_{ij})_{n \times p}$ and normalizing these order statistics, we get the standardized ranking matrix, 
\begin{align*}
	\R=
	\begin{pmatrix}
		\R_1\trans\\
		\vdots\\
		\R_n\trans
	\end{pmatrix} \defby  \left(\sqrt{\frac{12}{n^2-1}}(r_{ij}-\frac{n+1}{2}) \right)_{n \times p}.
\end{align*}
With the ranking matrix, \citet{spearman1904proof} introduced the famous Spearman's rank correlation matrix
\begin{align}
	\brho_n=\frac{1}{n} \R \trans \R =\frac{1}{n} \sum_{i=1}^n \R_i \R_i \trans.
\end{align}
From the statistical point of view, inducing the order statistics makes the correlation more robust to heavy-tailed distributions.  From the technical point of view, the ranking violates the independent structure among the rows of the observations. In particular, we have the following proposition.
\begin{proposition} \label{prob1}
	Assuming $\X_1, \ldots, \X_n,~i.i.d. \sim N(0,\bI_p)$, we have $\E \{\R \trans \R/n\} = \bI_p$ and $ \E \{\R  \R \trans/p\}=(n \bI_n-\one_n\one_n\trans)/(n-1)$. 
\end{proposition}

From Proposition \ref{prob1}, we can see that when the raw data is with i.i.d. rows, the ranking matrix $\R$ does not have independent rows anymore. Since the columns of the data matrix $\X$ is also i.i.d.,  the ranking matrix $\R$  has i.i.d. columns and in special, each column $(r_{1j},\ldots,r_{nj})$ is uniformly distributed from the permutations of $\{1,\ldots, n\}$.  \citet{bai2008large} studied this case and they proved that  the LSD of $\brho_n$ is also the Mar\u{c}enko-Pastur law. 

On the other hand, although we assume $\X_1, \ldots, \X_n,~i.i.d. \sim N(0,\bI_p)$, it actually covers a large class of distributions since the Spearman's correlation is rank-based. For example, when all the features are independent and continuous, we can transform each feature into a standard normal distribution where the ranking matrix is invariant.  It is referred to the monotonic invariance of rank-based correlation matrices introduced by \citet{weihs2018symmetric} and the non-paranormal model proposed by \citet{liu2009nonparanormal}. Thus, the independent case of previous works on Spearman's correlation matrix \cite[e.g.,][]{bai2008large,bao2015Spearman,bao2019tracy_sp} can be formulated as $\X_1, \ldots, \X_n,~i.i.d.~\sim N(0,\bI_p)$.

Assuming that the data matrix is with i.i.d. entries is very limited. A natural extension is to consider general covariance structures. For the multivariate normal distribution with general $\bSig$, we have the following result. 
\begin{proposition} \label{prob2}
	Assuming $\X_1, \ldots,\X_n,~i.i.d.\sim N(\mathbf{0},\bSig)$, where $\Sigma_{ii}=1,~i\in\{1,\ldots, p\}$, we have
	\begin{align*}
		\E \{\frac{1}{n} \R \trans \R\}=\frac{3}{n+1}\bSig_1+\frac{3(n-2)}{n+1} \bSig_2,~\E \{\frac{1}{p} \R  \R \trans\}=\frac{n}{n-1}(\bI_n-\frac{1}{n}\one_n\one_n\trans),
	\end{align*} 
	where 
	\begin{align*}
		\bSig_1=\left( \frac{2}{\pi} \arcsin(\Sigma_{ij}) \right)_{p \times p} \defby \frac{2}{\pi}\arcsin(\bSig),~\bSig_2=\left( \frac{2}{\pi} \arcsin(\Sigma_{ij}/2) \right)_{p \times p}\defby \frac{2}{\pi}\arcsin(\bSig/2).
	\end{align*}
\end{proposition}

From Proposition \ref{prob2}, we can see that the rows and columns of the ranking matrix $\R$ are both dependent if the data has a general covariance structure $\bSig$. From the perspective of random matrix theory, it is challenging to analyze such matrices.

\section{U-statistic of Spearman's correlation} \label{sec3}
For each feature $j \in \{1,\ldots,p\}$, we invoke the empirical cumulative distribution function
\begin{align*}
	\hat{F}_j(x)=\frac{1}{n}\sum_{i=1}^n I(x_{ij}\leq x),
\end{align*}
and the order statistics are
\begin{align*}
	r_{ij}=n \cdot \hat{F}_j(x_{ij})=1+\sum_{k\neq i} I(x_{kj}\leq x_{ij}), ~i\in\{1,\ldots, n\},~j\in\{1,\ldots, p\}.
\end{align*}
Throughout the paper, we consider the continuous population distribution and assume that there is no ties in data matrix.

Defining the sign vector
\begin{align*}
	\A_{ik}=\sign\left(\X_{i}-\X_{k}\right)=\left(\sign\left(x_{i1}-x_{k1}\right), \ldots, \sign\left(x_{ip}-x_{kp}\right)\right)\trans,
\end{align*}
we have
\begin{align*}
	r_{ij}-\frac{n+1}{2}=1+\sum_{k \neq i} \frac{1+ \sign\left(x_{ij}-x_{kj}\right)}{2}-\frac{n+1}{2}=\frac{1}{2}\sum_{k \neq i} \sign\left(x_{ij}-x_{kj}\right).
\end{align*}
In vector form, we can get
\begin{align} \label{Ri}
	\R_i=&	\sqrt{\frac{12}{n^2-1}}  \begin{pmatrix}
		r_{i1}-\frac{n+1}{2}\\
		\vdots\\
		r_{ip}-\frac{n+1}{2}
	\end{pmatrix}
	=\sqrt{\frac{3}{n^2-1}}\sum_{k\neq i}\A_{ik}. 
\end{align}
With these notations, we can rewrite the Spearman's rank correlation matrix as follows
\begin{align}
	\mathbf{\brho}_{n}=\frac{3}{n(n^2-1)} \sum_{i=1}^n \left\{ \sum_{k_1,k_2 \neq i}  \A_{i k_1} \A_{i k_2}\trans \right\}.
\end{align}
Following the definition of U-statistics, we can divide the summation into two parts and each part is a U-statistic. 

\begin{proposition} \label{prop3}
	Assuming the observations $\X_1,\ldots,\X_n$ have no ties, the Spearman's correlation matrix can be defined by	
	\begin{align}
		\brho_{n}=\frac{3}{n(n^2-1)} \sum_{i,j}^* \A_{ij} \A_{ij}\trans+\frac{3}{n(n^2-1)} \sum_{i,j,k}^* \A_{ij} \A_{ik}\trans,
	\end{align}
	where $\sum^*$ denotes summation over mutually different indices. 
\end{proposition}
When $\X_1,\ldots, \X_n$ are i.i.d. from a population with absolutely continuous densities, by Proposition \ref{prop3}, we can see 
\begin{align*}
	\E (\brho_{n})=&\frac{3}{n(n^2-1)} n(n-1) \cdot \cov\left(\A_{12}\right)+\frac{3}{n(n^2-1)} n(n-1)(n-2)\cdot\cov\left(\A_1\right)
	=\frac{3}{n+1}  \cov\left(\A_{12}\right)+\frac{3(n-2)}{n+1} \cov\left(\A_1\right),
\end{align*}
where
\begin{align}\label{ai}
	\A_i  \defby \E \{\A_{ik}\mid \X_i\}.
\end{align}
The population Spearman's correlation matrix consists of two covariance matrices $\cov\left(\A_{12}\right)$ and $\cov\left(\A_1\right)$.  For the classical sample covariance matrix
\begin{align*}
	\bS_n=\frac{1}{2n(n-1)}  \sum_{i,j}^* (\X_i-\X_j)(\X_i-\X_j)\trans,
\end{align*}
we know $\cov(\X_i-\X_j)=2 \cov(\X_i)$. Thus, the $\sign$ function introduces a non-linear correlation into rank-based correlation matrices and these two covariance matrices are different, e.g., 
\begin{align*}
	\cov\left(\A_{12}\right)=\bSig_1=\frac{2}{\pi}\arcsin(\bSig),~\mbox{and}~\cov\left(\A_1\right)=\bSig_2=\frac{2}{\pi}\arcsin(\bSig/2),
\end{align*}
for the multivariate normal distribution $N(\mathbf{0},\bSig)$.

Since both $\cov\left(\A_{12}\right)$ and $\cov\left(\A_1\right)$ can describe correlations among features, 
an intuitive way is to conduct U-statistics for these two covariance matrices separately, e.g., 
\begin{align} \label{kendall}
	\btau_n=\frac{1}{n(n-1)}  \sum_{i,j}^* \A_{ij} \A_{ij}\trans,
\end{align} 
and
\begin{align} \label{newsp}
	\tilde{\brho}_n=\frac{3}{n(n-1)(n-2)} \sum_{i,j,k}^* \A_{ij} \A_{ik}\trans.
\end{align}
Interestingly, $\btau_n$ is exactly Kendall's rank correlation matrix proposed by \cite{kendall1938new} and $\tilde{\brho}_n$ is an improved Spearman's rank correlation matrix proposed by  \cite{hoeffiding1948class}. See also Example 3 of \cite{han2017distribution}.  Our recent work \cite{li2021eigenvalues} studied the LSD of $\btau_n$ and in this work, we consider Spearman's correlation matrices $\brho_n$ and $\tilde{\brho}_n$. 

\section{Limiting spectral distribution} \label{sec4}
For an $n \times n$ Hermitian matrix $\bH_n$ whose eigenvalues are $\lambda_1,\ldots,\lambda_n$, the empirical spectral distribution of $\bH_n$ is defined as
\begin{align*}
	F^{\bH_n}(x)=\frac{1}{n}\sum_{i=1}^n I(\lambda_i\leq x).
\end{align*}
The limit of $F^{\bH_n}$ is called the limiting spectral distribution of $\bH_n$.  In random matrix theory, the LSD is usually defined by its Stieltjes transform
\begin{align*}
	s_{F}(z)=\int\frac{1}{x-z} dF(x),~z \in \mathbb{C}^{+},
\end{align*}
where $\mathbb{C}^{+}$ denotes the upper complex plane. With a Stieltjes transform $s(z)$, the distribution function can be  obtained by the inversion formula
\begin{align*}
	F(b)-F(a)=\lim_{\nu \to 0^+}\frac{1}{\pi}\int_a^b \Im s_{F}(x+i\nu)dx,
\end{align*}
where $\Im(\cdot)$ is the imaginary part and $i$ is the imaginary unit.

Note
\begin{align*} 
	\tilde{\brho}_n=\frac{3}{n} \sum_{i=1}^n \left\{\frac{1}{(n-1)(n-2)} \sum_{j,k \neq i}^* \A_{ij} \A_{ik}\trans  \right\}
\end{align*}
and for a given $i$, 
\begin{align*}
	\frac{1}{(n-1)(n-2)} \sum_{j,k \neq i}^* \A_{ij} \A_{ik}\trans 
\end{align*}
is also a U-statistic for $\A_i \A_i \trans$. Intuitively,  the improved Spearman's rank correlation matrix $\tilde{\brho}_n$ is close to the random matrix
\begin{align}
	\W_n=\frac{3}{n} \sum_{i=1}^n \A_i\A_i\trans.
\end{align}
The following result shows that $\tilde{\brho}_n$ and $\W_n$ share the same limiting spectral distribution.
\begin{theorem} \label{thm1}
	Assuming $\X_1, \ldots, \X_n,~i.i.d.\sim N(\mathbf{0},\bSig)$  where $\Sigma_{ii}=1,~i\in\{1,\ldots, p\}$, and if
	\begin{enumerate}[label=$\mathrm{(\roman*)}$]
		\item $p/n \to y \in (0, \infty)$;
		\item $\tr\left(\bSig^2 \right)=o(p^2)$;
\end{enumerate}
	we have
	\begin{align}
		L\left(F^{\tilde{\brho}_n}, F^{\W_n}\right)  \to 0,~\text{in probability},
	\end{align}
	where $L(\cdot,\cdot)$ is the Levy distance between two distribution functions. 
\end{theorem}

For Spearman's rank correlation matrix $\brho_n$,  by Proposition \ref{prop3}, we know
\begin{align*}
	\brho_n=\frac{3}{n+1} \btau_n+\frac{n-2}{n+1} \tilde{\brho}_n=\frac{3}{n+1} \btau_n-\frac{3}{n+1} \tilde{\brho}_n+\tilde{\brho}_n,
\end{align*}
where $\btau_n$ is the Kendall's rank correlation matrix \eqref{kendall}. We claim that $\brho_n$ and $\W_n$ also share the same limiting spectral distribution.   The following result establishes this claim.
\begin{proposition}[Weak convergence] \label{prop4}
	Under the assumptions of Theorem \ref{thm1}, we have
	\begin{align}
		L\left(F^{\brho_n}, F^{\W_n}\right)  \to 0,~\text{in probability}.
	\end{align}
\end{proposition}
The proofs of Theorem \ref{thm1} and Proposition \ref{prop4} are based on Corollary A.41 of \citet{bai2010spectral} and we need to show
$\E \left\|\tilde{\brho}_n-\W_n\right\|_2^2=o(p)$ and $\E \left\|\brho_n-\W_n\right\|_2^2=o(p)$. To further refine the conclusion, e.g., to show the strong convergence, we need to prove $ \left\|\tilde{\brho}_n-\W_n\right\|_2^2/p \to 0,a.s.$ which is challenging since the calculation of higher order moments of $\A_i$ is hard. Here, we turn to the ranking statistics $\R$. Intuitively, 
\begin{align*}
	\R=\left(\sqrt{\frac{12}{n^2-1}}(r_{ij}-\frac{n+1}{2}) \right)_{n \times p}= \sqrt{\frac{12 n^2 }{n^2-1}} \left(\hat{F}_j(x_{ij})-\frac{n+1}{2n} \right)_{n \times p}	\approx 2\sqrt{3} \left(F_j(x_{ij})-\frac{1}{2} \right)_{n \times p}
	=\sqrt{3} (\A_1,\ldots,\A_n)\trans.
\end{align*}
Thus, we can control the difference between $F^{\brho_n}$ and $F^{\W_n}$ by bounding $\R-\sqrt{3} (\A_1,\ldots,\A_n)\trans$. Specifically, we will use Corollary A.42 of \citet{bai2010spectral}, i.e., 
\begin{align*}
	L^4\left(F^{\brho_n}, F^{\W_n} \right) \leq & \frac{2 \tr \left(\brho_n+ \W_n \right) }{n p^2} \|\R-\sqrt{3} (\A_1,\ldots,\A_n)\trans\|_2^2
\end{align*}
and then we can show the strong convergence as follows. 
\begin{theorem}[Strong convergence]  \label{thm1-1}
	Assuming $\X_1, \ldots, \X_n$ are i.i.d. continuous random vectors and $p/n \to y \in (0, \infty)$, we have 
	\begin{align}
		L\left(F^{\brho_n}, F^{\W_n}\right)  \to 0,~\text{almost surely}.
	\end{align}
\end{theorem}
It is noted that Theorem \ref{thm1-1} provides a stronger conclusion with a weaker condition and the key technical tool is the Dvoretzky–Kiefer–Wolfowitz inequality for empirical cumulative distribution functions. By Proposition \ref{prop4} or Theorem \ref{thm1-1},  to study the LSD of $\brho_n$, we can consider the matrix $\W_n$ which is the sample covariance matrix of the random vectors $\sqrt{3}(\A_1,\ldots,\A_n)$. It is noted that $\A_1,\ldots, \A_n$ are i.i.d., and the LSD can be derived from classical results on sample covariance matrices, e.g., Theorem 1 of \citet{bai2008large}.  	Assuming $\X_1, \ldots, \X_n,~i.i.d. \sim N(\mathbf{0},\bSig)$ where $\Sigma_{ii}=1$, we have
\begin{align*}
	\A_i=\E\{\A_{ij}\mid \X_i\}=\begin{pmatrix}
		2 \Phi(x_{i1})-1\\
		\vdots\\
		2 \Phi(x_{ip})-1\\
	\end{pmatrix}.
\end{align*}
Here $\Phi(\cdot)$ is the cumulative distribution function of $N(0,1)$. Therefore, each entry of $(\A_1,\ldots,\A_n)$ follows the uniform distribution $U[-1,1]$ and 
\begin{align*}
	\cov(\A_i)=\bSig_{2}=\frac{2}{\pi}\arcsin(\bSig/2).
\end{align*}
Based on the main result of \citet{bai2008large}, we can derive the LSD of Spearman's rank correlation matrix as follows.
\begin{theorem}\label{thm2}
	Assume $\X_1, \ldots, \X_n,~i.i.d. \sim N(\mathbf{0},\bSig)$ where $\Sigma_{ii}=1, ~i \in \{1,\ldots, p\}$, and 
	\begin{enumerate}[label=$\mathrm{(\roman*)}$]
		\item $p/n \to y \in (0, \infty)$ as $n \to \infty$;
		\item the spectral norm of $\bSig$ is uniformly bounded by a constant $C$;
		\item the empirical spectral distribution of $ 6 \arcsin(\bSig/2)/\pi$ tends to a non-random probability distribution $H$.  
	\end{enumerate}
	Then, with probability 1, $F^{\brho_n}$ tends to a probability distribution, whose Stieltjes transform $m=m(z),~z \in \mathbb{C}^{+}$ satisfies 
	\begin{align}
		m=\int \frac{1}{t(1-y-yzm)-z} dH(t).
	\end{align}
\end{theorem}

Noting the Spearman's correlation matrix is rank-based, Theorem \ref{thm2} actually holds for the non-paranormal distribution proposed by \citet{liu2009nonparanormal}. In particular, for a random vector $\Z=(z_1,\ldots,z_p)\trans \in \mR^p$, there exist monotone functions $\left\{f_j\right\}_{j=1}^p$ such that $\left( f_1(z_1),\ldots,f_p(z_p)\right) \sim N(\mathbf{0},\bSig)$ where $\Sigma_{ii}=1, ~i \in \{1,\ldots, p\}$. Denoting the cumulative distribution function of $z_j$ as $F_j$, \citet{liu2009nonparanormal} show that $f_j(t)=\Phi^{-1}(F_j(t))$ and the non-paranormal distribution model actually assumes
	\begin{align*}
		\left( \Phi^{-1}(F_1(z_1)),\ldots,\Phi^{-1}(F_p(z_p))\right) \sim N(\mathbf{0},\bSig).
	\end{align*}
We can see that the non-paranormal distribution extends the multivariate normal distribution and more discussions can be found in \citet{lu2021robust}.  More generally, we can further extend the assumption to general distributions with some moment conditions and it is referred to  \citet{li2021eigenvalues} for more details. To illustrate the robustness of Spearman's rank correlation, we consider a toy example where each element follows the Cauchy distribution.
\begin{example}
		For the data matrix
	\begin{align*}
		\X=
		\begin{pmatrix}
			\X_1\trans\\
			\vdots\\
			\X_n\trans
		\end{pmatrix}=\begin{pmatrix}
			x_{11} & \cdots & x_{1p}\\
			\vdots &\ddots& \vdots\\
			x_{n1} &\cdots & x_{np}
		\end{pmatrix}=
		(x_{ij})_{n \times p},
	\end{align*}
where $x_{ij}$ follows the standard Cauchy distribution, we consider the monotone transformation 
\begin{align*}
	x_{ij} \mapsto \Phi^{-1}\left(\frac{1}{2}+\frac{1}{\pi} \arctan(x_{ij})\right) \defby y_{ij}. 
\end{align*}
As we know, the transformation does not change the ranking matrix $\R$ and to analyze $\R$, we can make assumptions on $y_{ij}$ without loss of generality.  For the independent case,  i.e., $x_{i1},\ldots, x_{ip}$ are independent,  we have $(y_{i1},\ldots, y_{ip})\trans \sim N(0,\bI)$. For the dependent case, we need further assumptions to define the correlation structure. Noting $y_{ij} \sim N(0,1)$, a natural way is to assume
\begin{align}\label{exfor1}
	(y_{i1},\ldots, y_{ip})\sim N(\mathbf{0},\bSig) 
\end{align}
where $\Sigma_{ii}=1, ~i \in \{1,\ldots, p\}$. This is exactly the non-paranormal distribution and $\bSig$ is called the latent generalized correlation matrix by \citet{lu2021robust}. In summary, we have the relationships   
\begin{align*} 
	\text{Spearman's correlation}~(x_{ij},x_{ik})=	\text{Spearman's correlation}~ (y_{ij},y_{ik})=\frac{2}{\pi} \arcsin\left(\text{Pearson's correlation}~(y_{ij},y_{ik}) \right),
\end{align*}
where the second equality is due to the non-paranormal distribution assumption \eqref{exfor1}.
\end{example}
\section{Discussions and simulations} \label{sec5}
In this section, we compare several important covariance/correlation matrices.   Given the data $\X_1,\ldots, \X_n \in \mR^p$, there are four classical covariance/correlation matrices in statistical applications:
\begin{itemize}
	\item Sample covariance matrix:
	\begin{align*}
		\bS_n=\frac{1}{n-1}\sum_{i=1}^n (\X_i-\bar{\X})(\X_i-\bar{\X})\trans=\frac{1}{2n(n-1)}  \sum_{i,j}^* (\X_i-\X_j)(\X_i-\X_j)\trans= \frac{1}{n-1}\sum_{i=1}^n \X_i \X_i\trans-\frac{n}{n-1} \bar{\X} \bar{\X}\trans,
	\end{align*}
where $\bar{\X}=\frac{1}{n} \sum_{i=1}^n \X_i$;
	\item Pearson's correlation matrix \cite{pearson1909determination}:  
	\begin{align*}
		\mathbf{P}_n=\{\diag(\bS_n)\}^{-1/2} \bS_n \{\diag(\bS_n)\}^{-1/2},
	\end{align*}
	where $\diag(\cdot)$ is a diagonal matrix with the diagonal entries of the matrix;
	\item Kendall's correlation matrix \cite{kendall1938new}:
	\begin{align*}
		\btau_n=\frac{1}{n(n-1)}  \sum_{i,j}^* \A_{ij} \A_{ij}\trans,
	\end{align*}
	where $\A_{ij}=\sign(\X_i-\X_j)$;
	\item Spearman's correlation matrix \cite{spearman1904proof}
	\begin{align*}
		\brho_{n}=\frac{3}{n(n^2-1)} \sum_{i,j}^* \A_{ij} \A_{ij}\trans+\frac{3}{n(n^2-1)} \sum_{i,j,k}^* \A_{ij} \A_{ik}\trans,
	\end{align*}
	and the improved Spearman's correlation matrix \cite{hoeffiding1948class}
	\begin{align*} 
		\tilde{\brho}_n=\frac{3}{n(n-1)(n-2)} \sum_{i,j,k}^* \A_{ij} \A_{ik}\trans.
	\end{align*}
\end{itemize}

For comparison purposes, we assume that $\X_1, \ldots, \X_n,~i.i.d. \sim N(\mathbf{0},\bSig)$. Since $\mathbf{P}_n, \btau_n$ and $\brho_n, \tilde{\brho}_n$ are all correlation matrices which are invariant to the scale of features, we assume $\Sigma_{ii}=1$ to make a fair comparison on the sample covariance matrix.  

\begin{figure}[!ht]
	\centerline{
		\begin{tabular}{cccc}	
			\psfig{figure=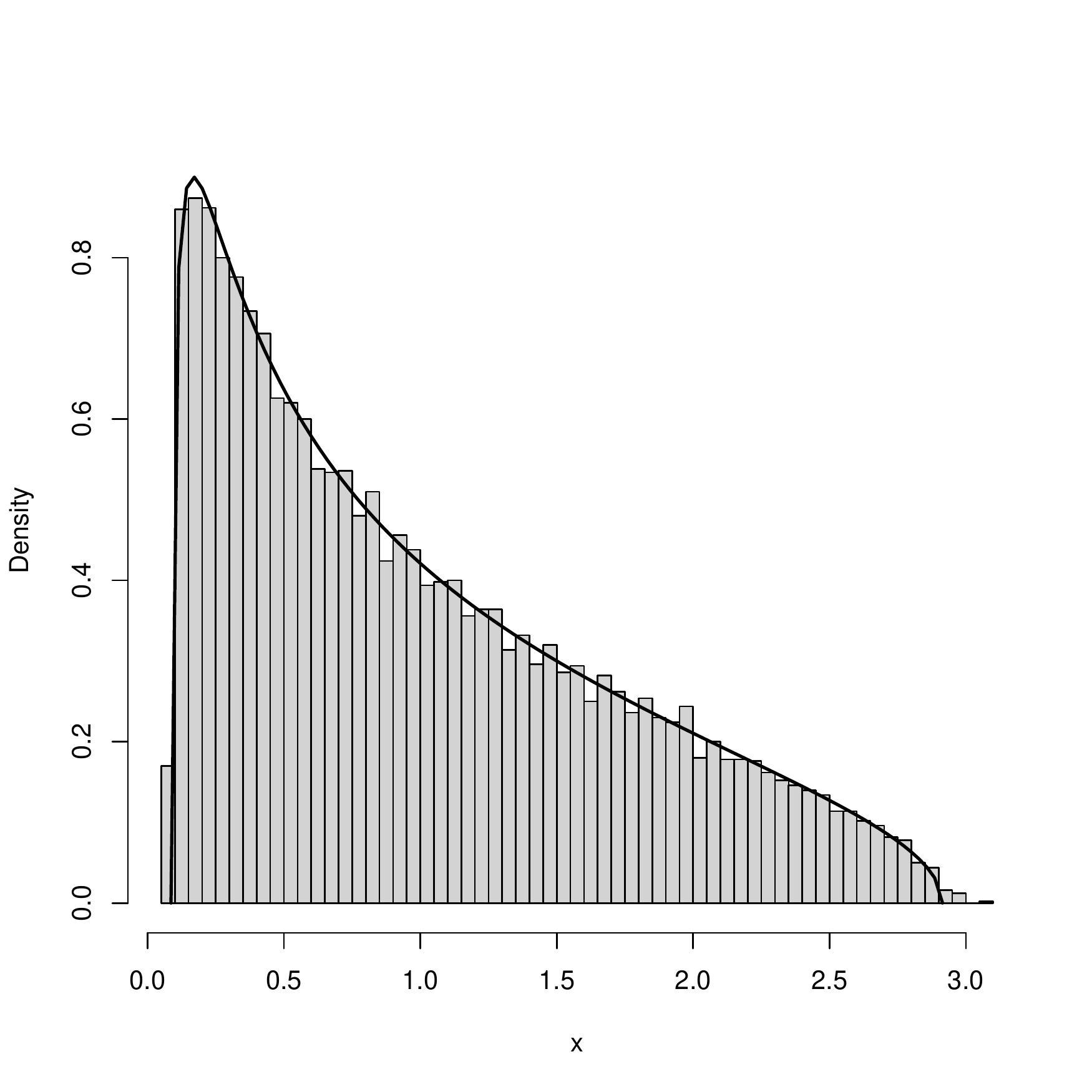,width=1.5 in,angle=0} &
			\psfig{figure=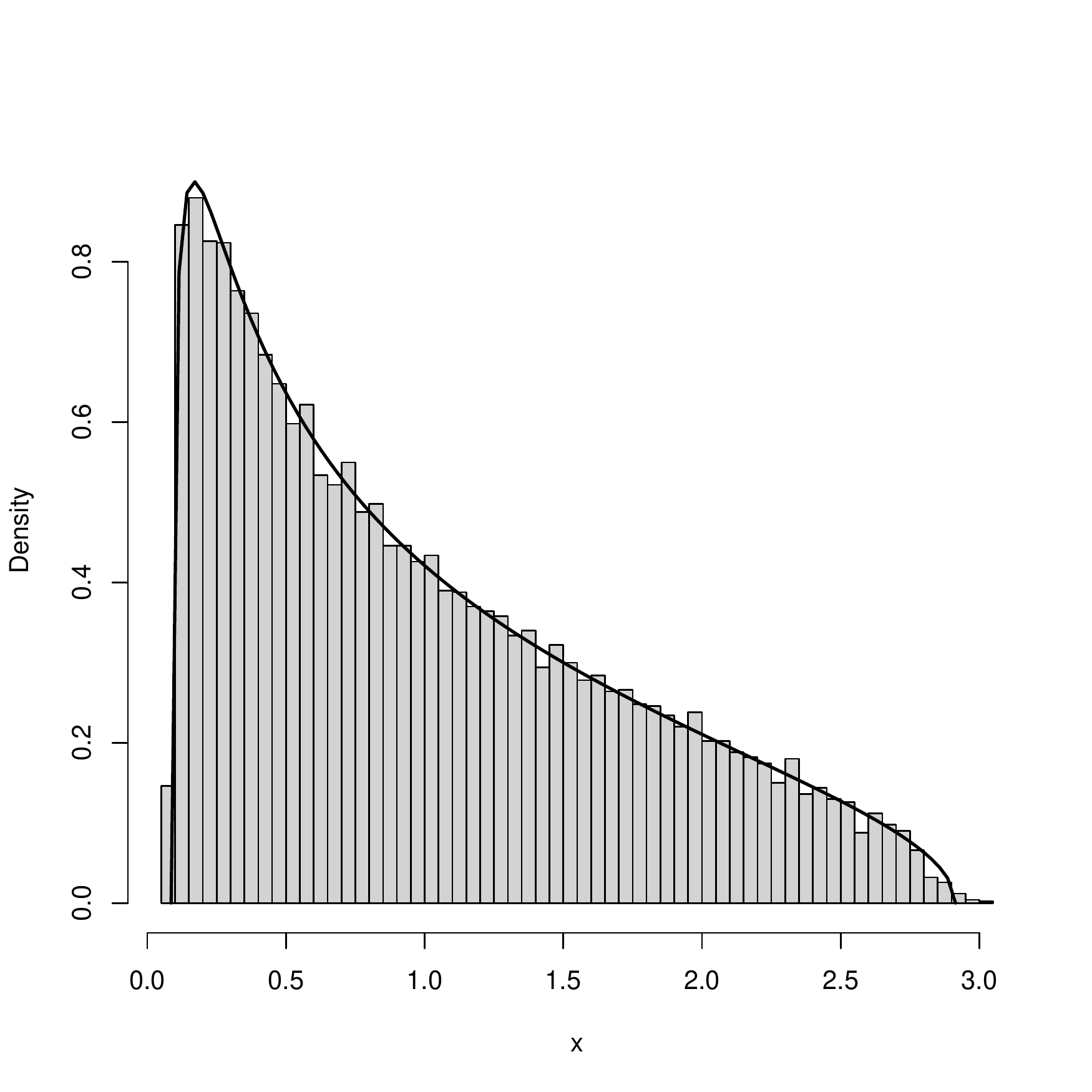,width=1.5 in,angle=0} & 
			\psfig{figure=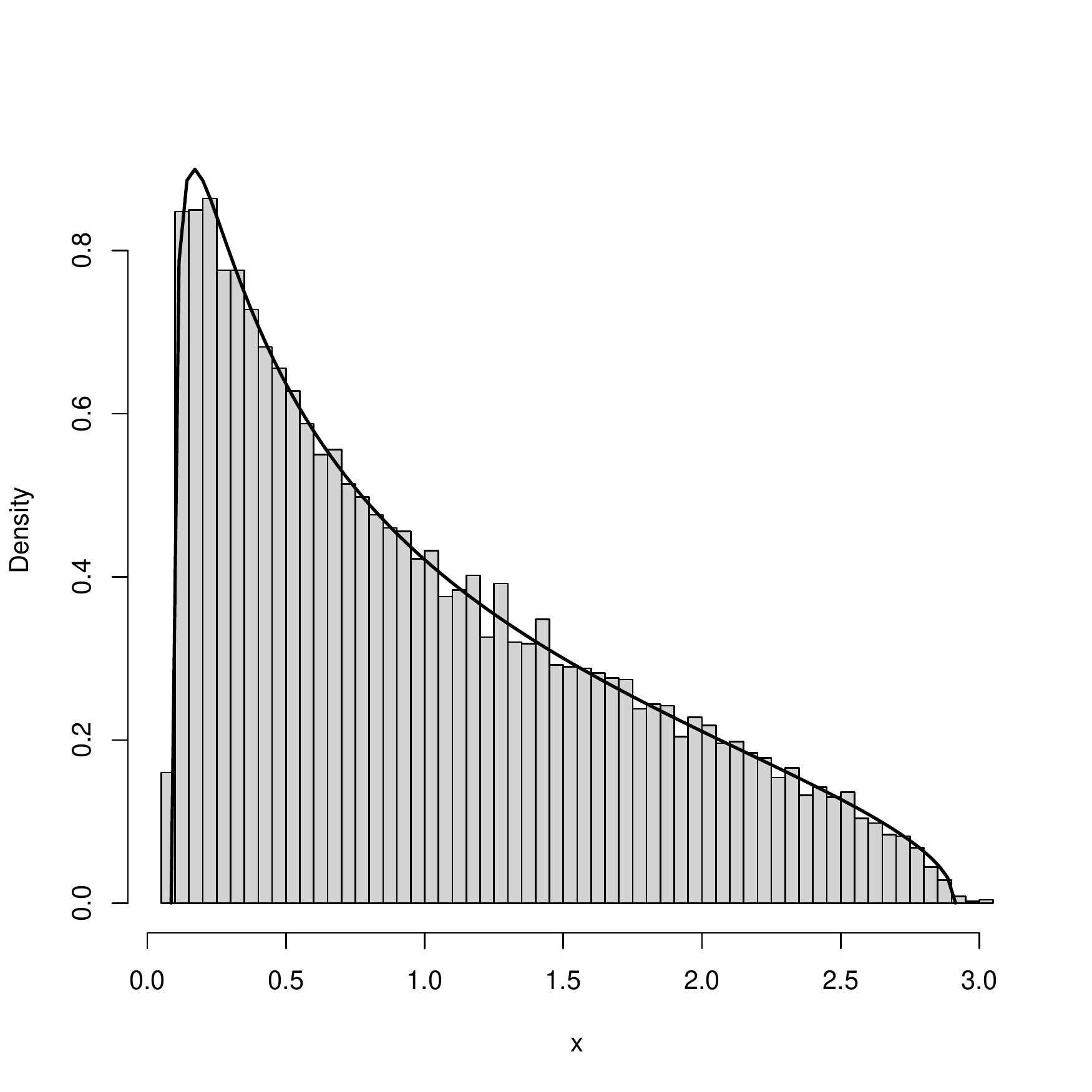,width=1.5 in,angle=0} &
			\psfig{figure=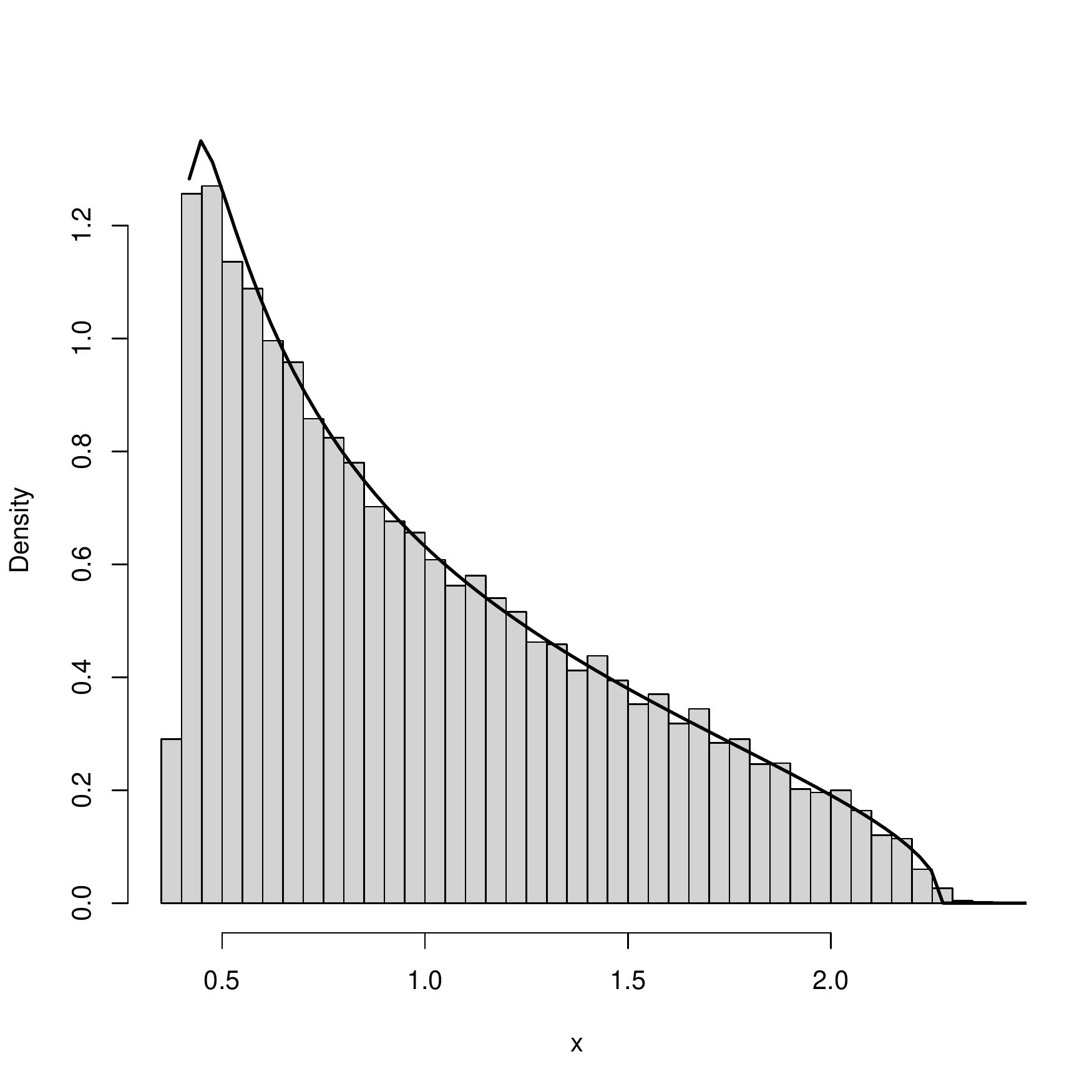,width=1.5 in,angle=0} \\			
			\multicolumn{4}{c}{$(n,p)=(200,100)$}\\	
			\psfig{figure=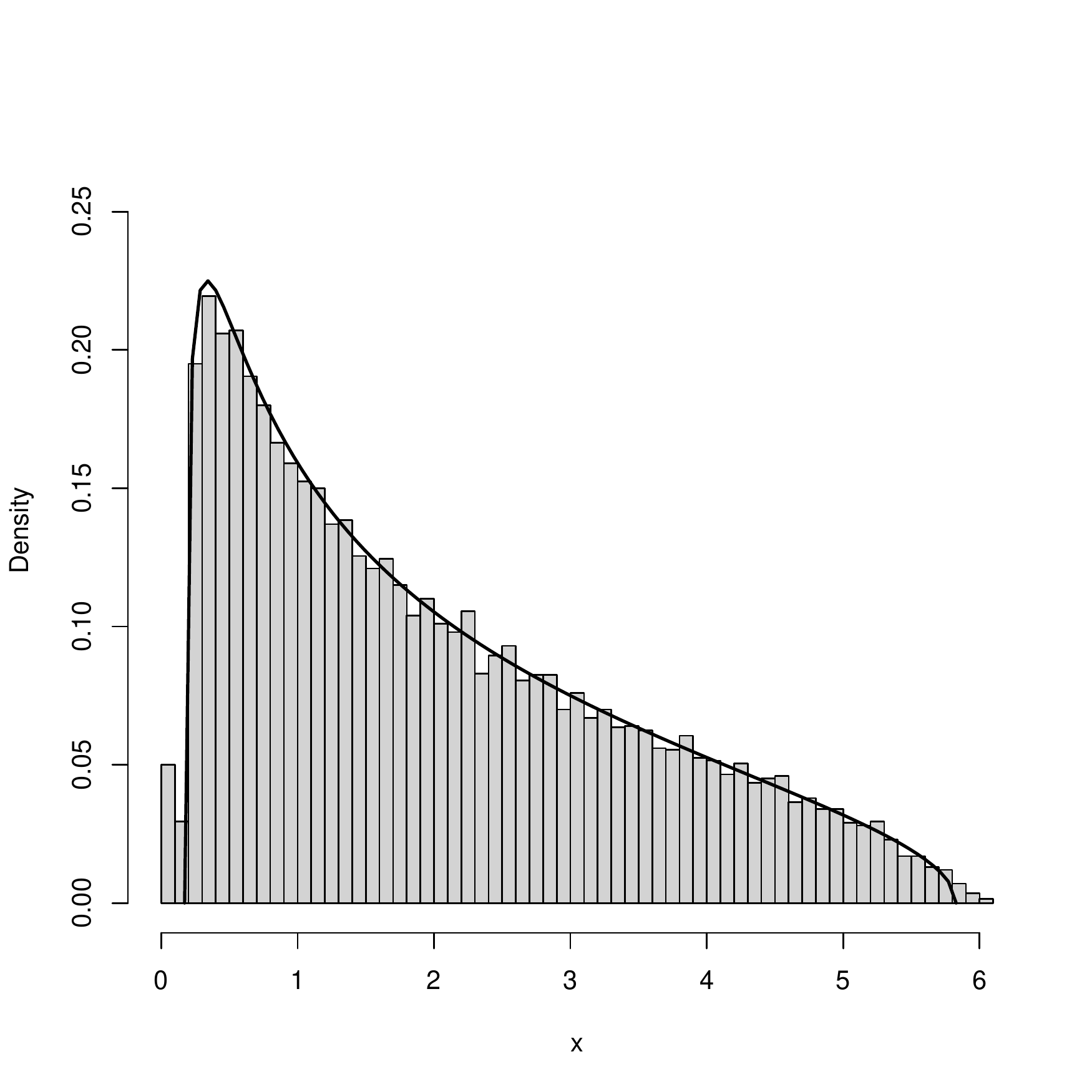,width=1.5 in,angle=0} &
			\psfig{figure=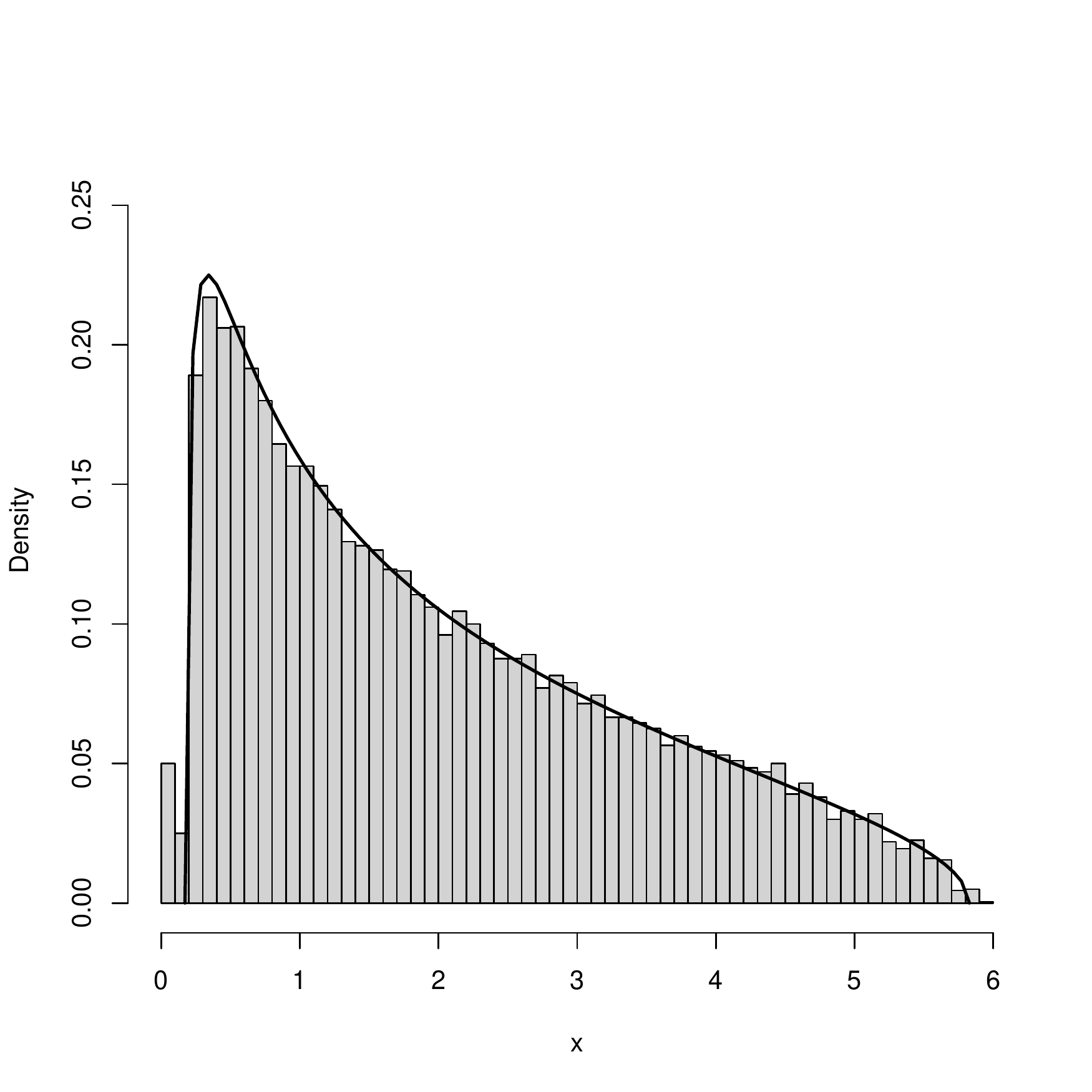,width=1.5 in,angle=0} & 
			\psfig{figure=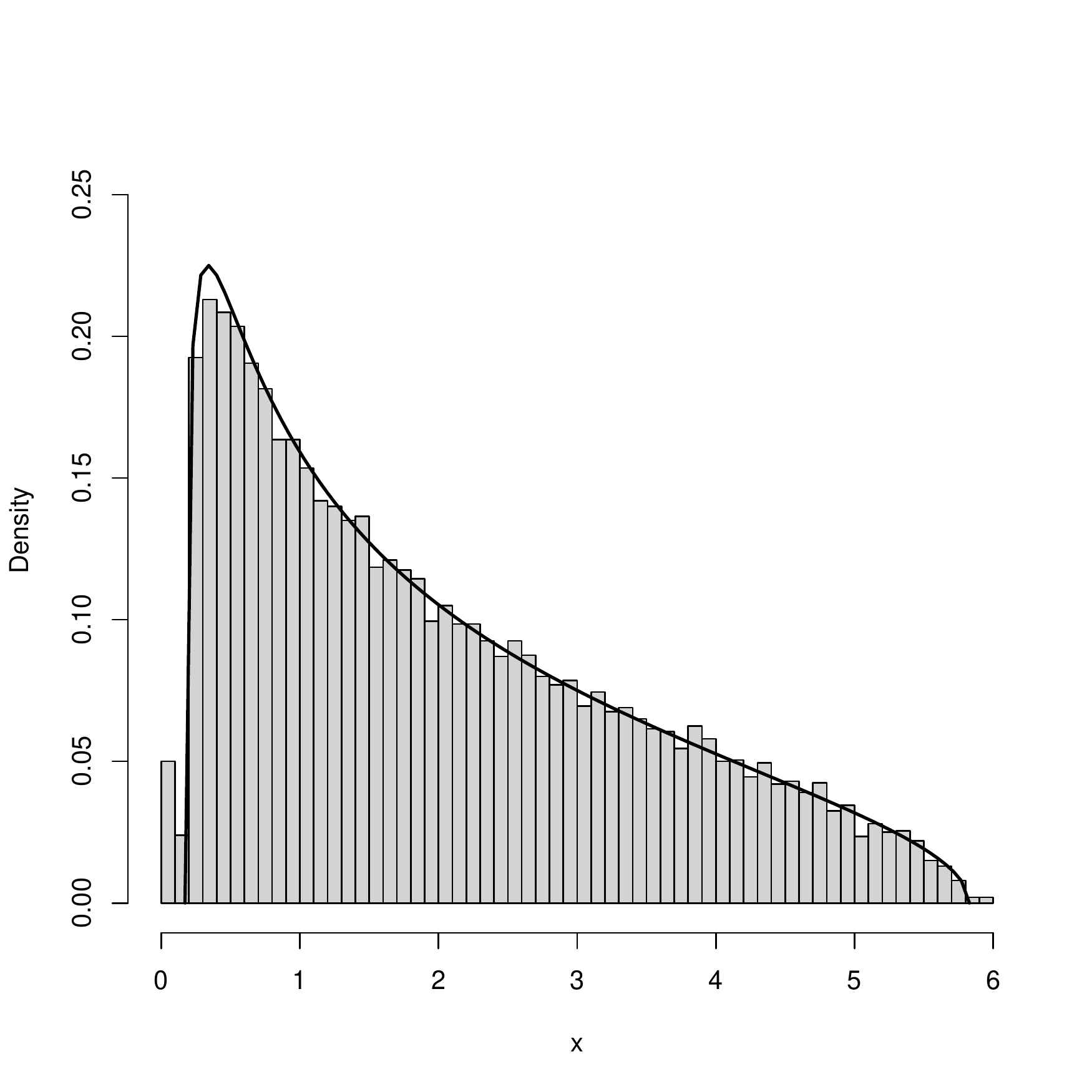,width=1.5 in,angle=0} &
			\psfig{figure=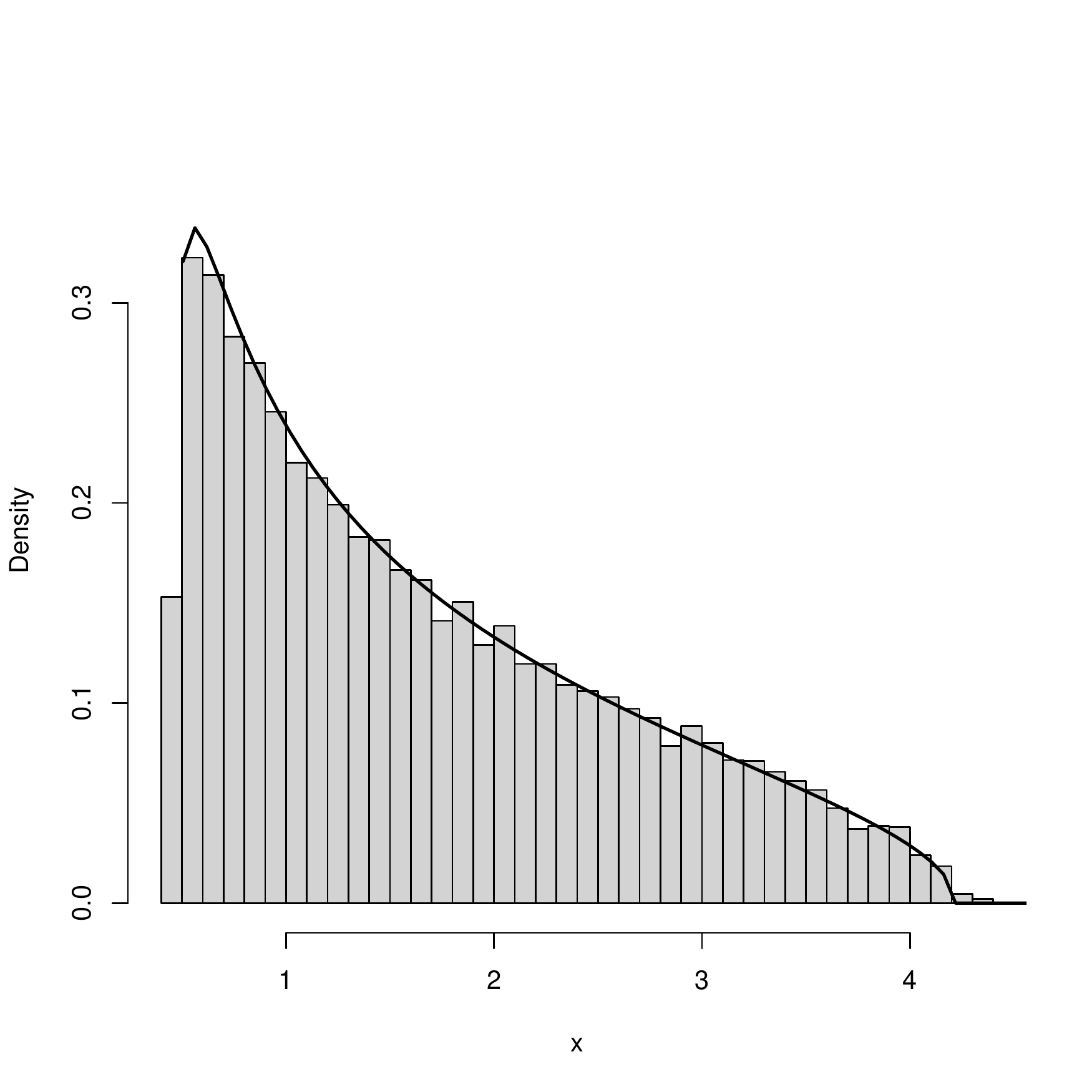,width=1.5 in,angle=0} \\
			\multicolumn{4}{c}{$(n,p)=(100,200)$}\\	
			Covariance& Pearson& Spearman& Kendall\\
		\end{tabular}
	}
	\caption{The limiting spectral distributions of sample covariance matrix, Pearson's correlation matrix, Spearman's correlation matrix and Kendall's correlation matrix where  $\X_1, \ldots, \X_n,~i.i.d. \sim N(0,\bI)$.}
	\label{fig1}
\end{figure}

When $\bSig=\bI$, \citet{marvcenko1967distribution} showed that the LSD of $\bS_n$ is the standard Mar\u{c}enko-Pastur law which has the density function 
\begin{align*}
	p(x)=\frac{1}{2\pi x y} \sqrt{(x_{+}-x) (x-x_{-})} I(x_{-} \leq x \leq x_{+}),
\end{align*}
for $y \leq 1$ and has a point mass $1-1/y$ at the origin if $y>1$. Here $p/n \to y \in (0, \infty)$, $x_{+}=(1+\sqrt{y})^2$ and $x_{-}=(1-\sqrt{y})^2$. \citet{jiang2004limiting} proved that the LSD of $\mathbf{P}_n$ is also the Mar\u{c}enko-Pastur law.  \citet{bai2008large} established the Mar\u{c}enko-Pastur law for the Spearman's rank correlation matrix. \citet{bandeira2017marvcenko} proved that the LSD of Kendall's correlation matrix is an affine transformation of the Mar\u{c}enko-Pastur law, i.e., 
\begin{align*}
	\frac{2}{3} Y+\frac{1}{3},
\end{align*}
where $Y$ follows the Mar\u{c}enko-Pastur law. Fig \ref{fig1} shows the theoretical LSDs and the empirical distributions based on 100 replications for these four matrices.

\begin{figure}[!ht]
	\centerline{
		\begin{tabular}{cccc}	
			\psfig{figure=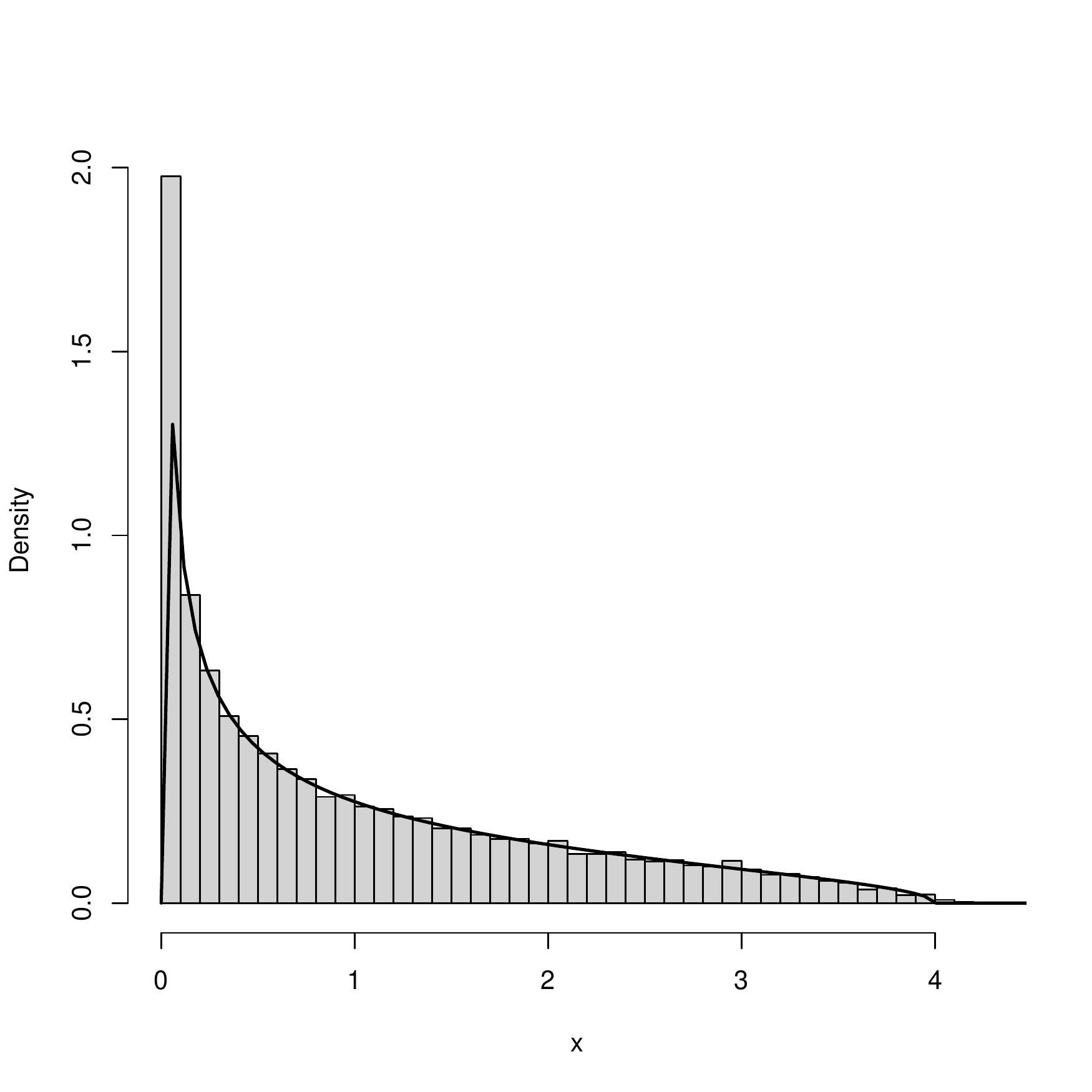,width=1.5 in,angle=0} &
			\psfig{figure=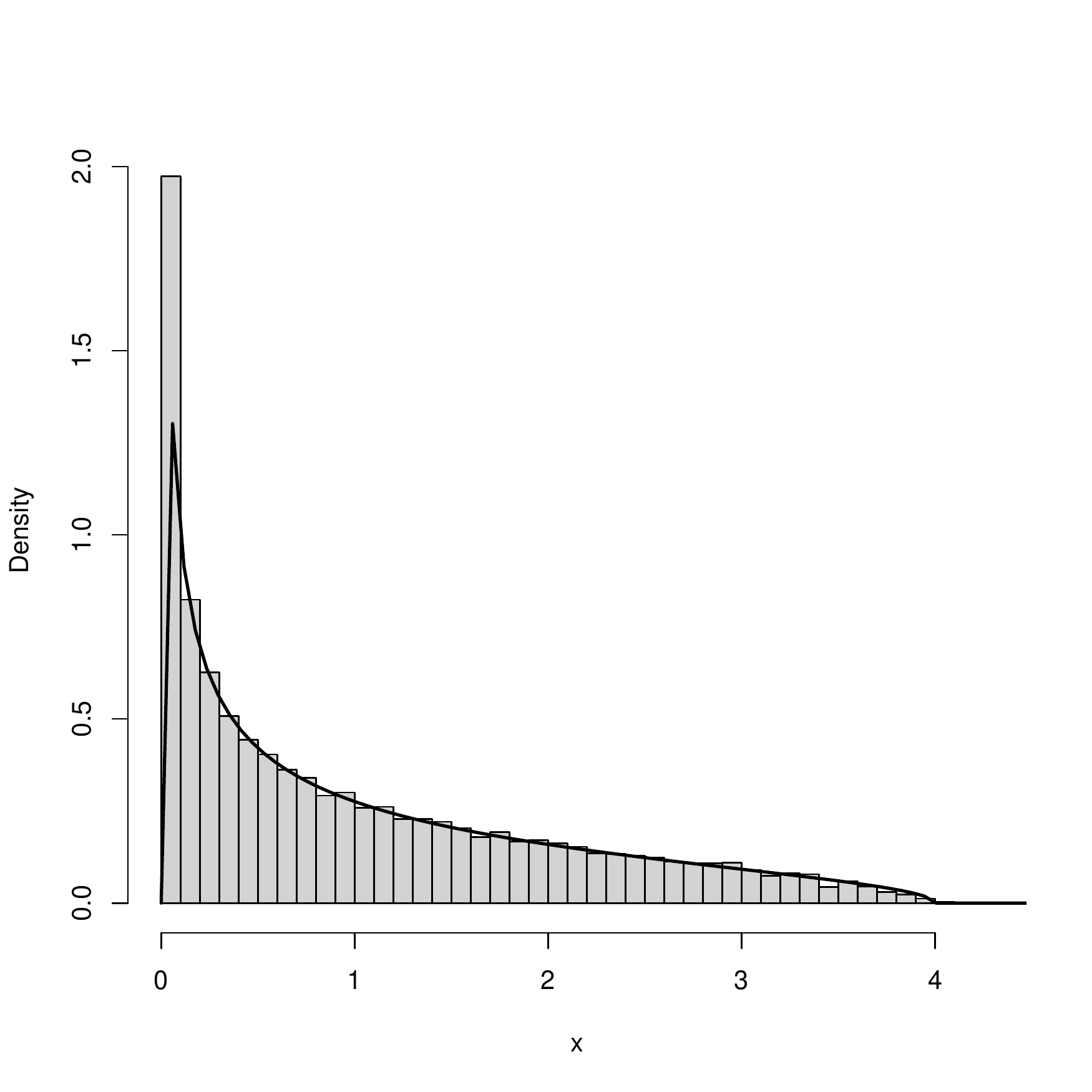,width=1.5 in,angle=0} & 
			\psfig{figure=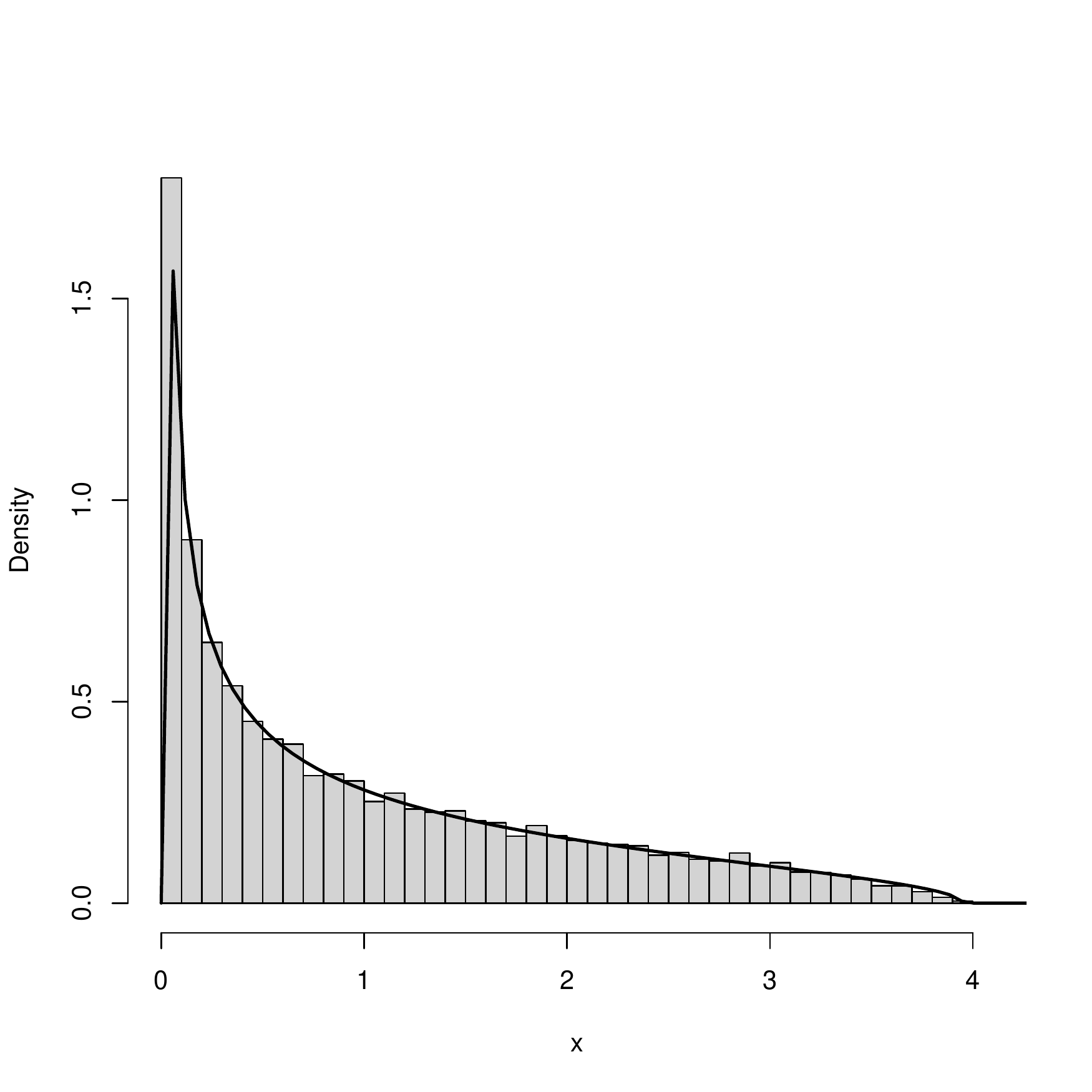,width=1.5 in,angle=0} &
			\psfig{figure=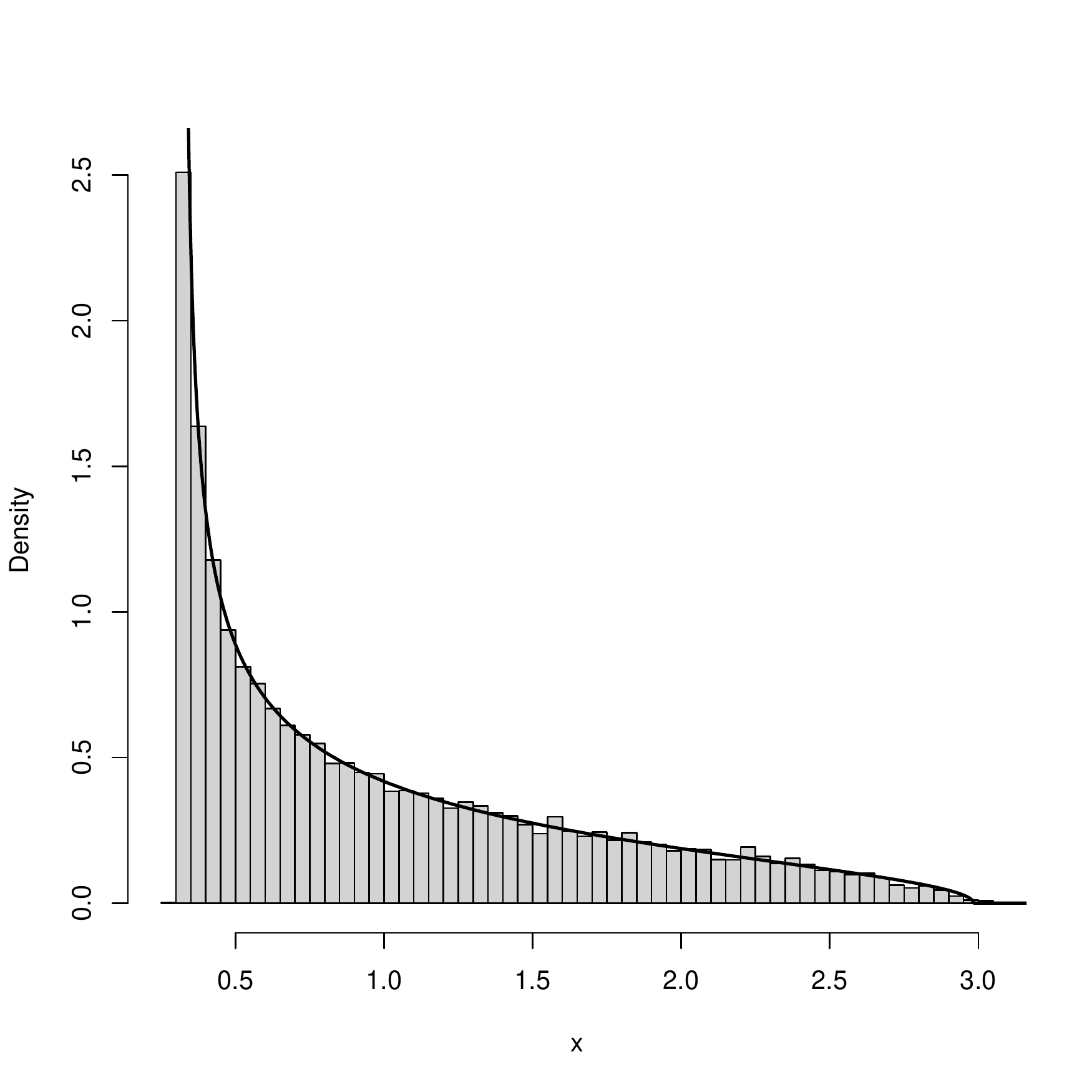,width=1.5 in,angle=0} \\			
			\multicolumn{4}{c}{$(n,p)=(200,100)$}\\	
			\psfig{figure=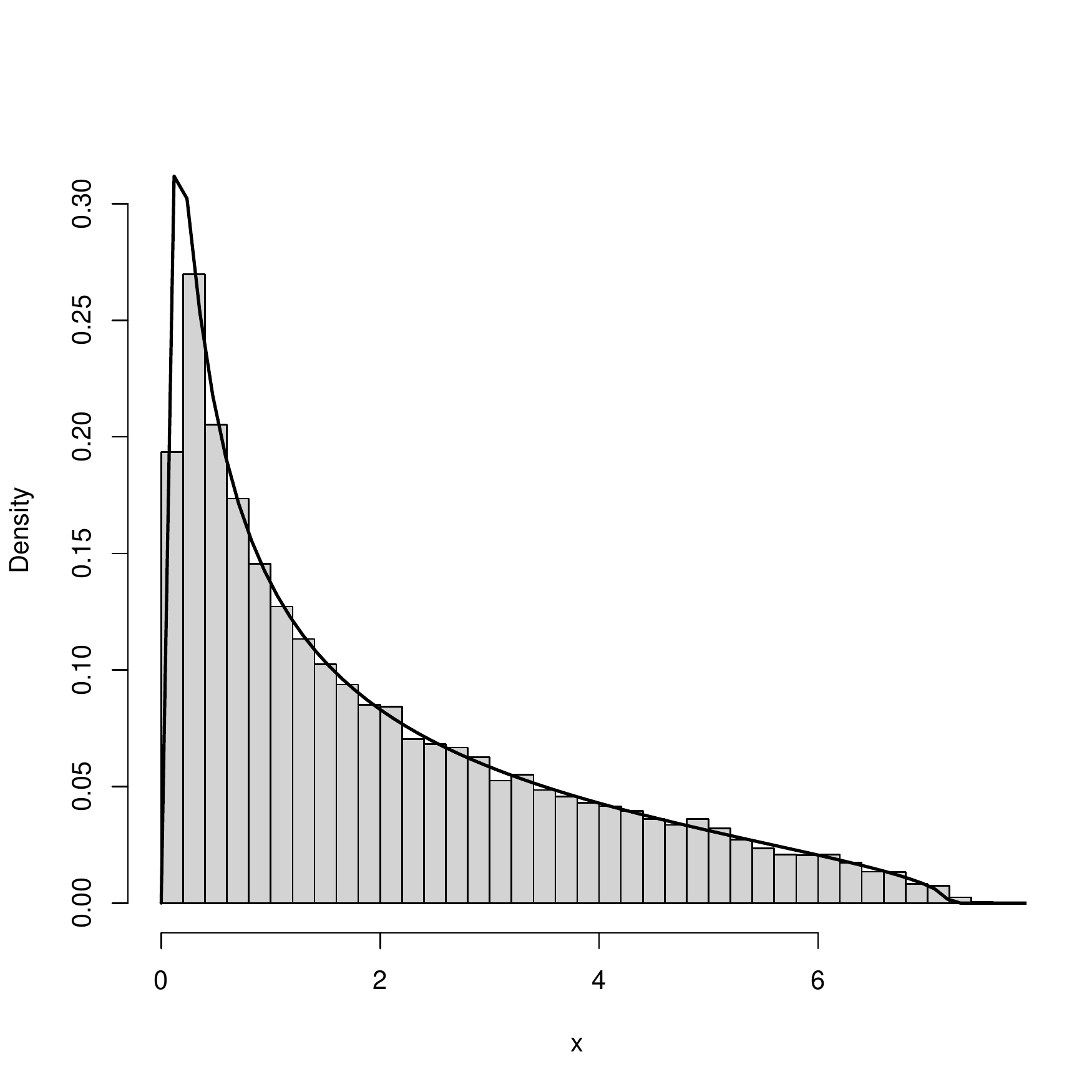,width=1.5 in,angle=0} &
			\psfig{figure=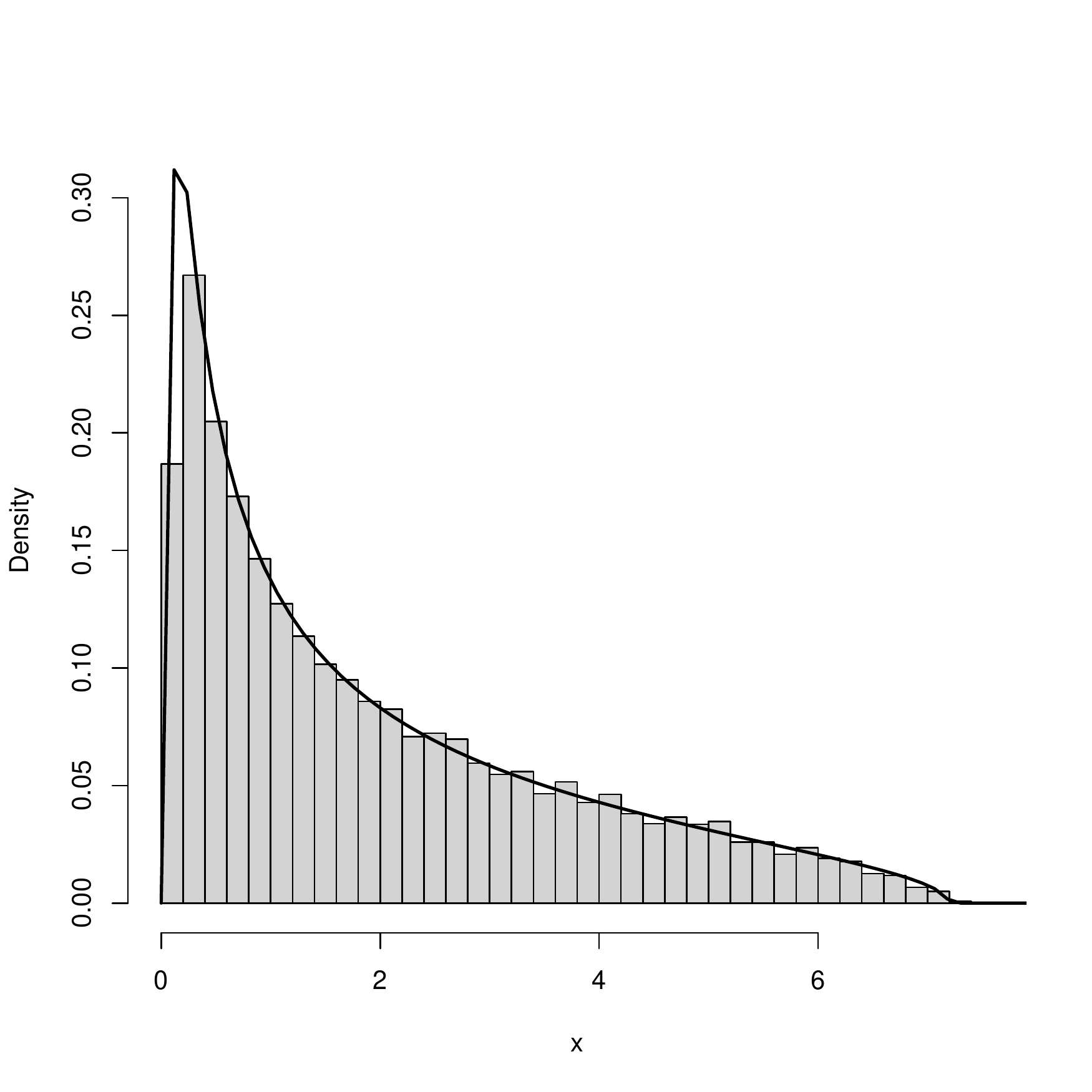,width=1.5 in,angle=0} & 
			\psfig{figure=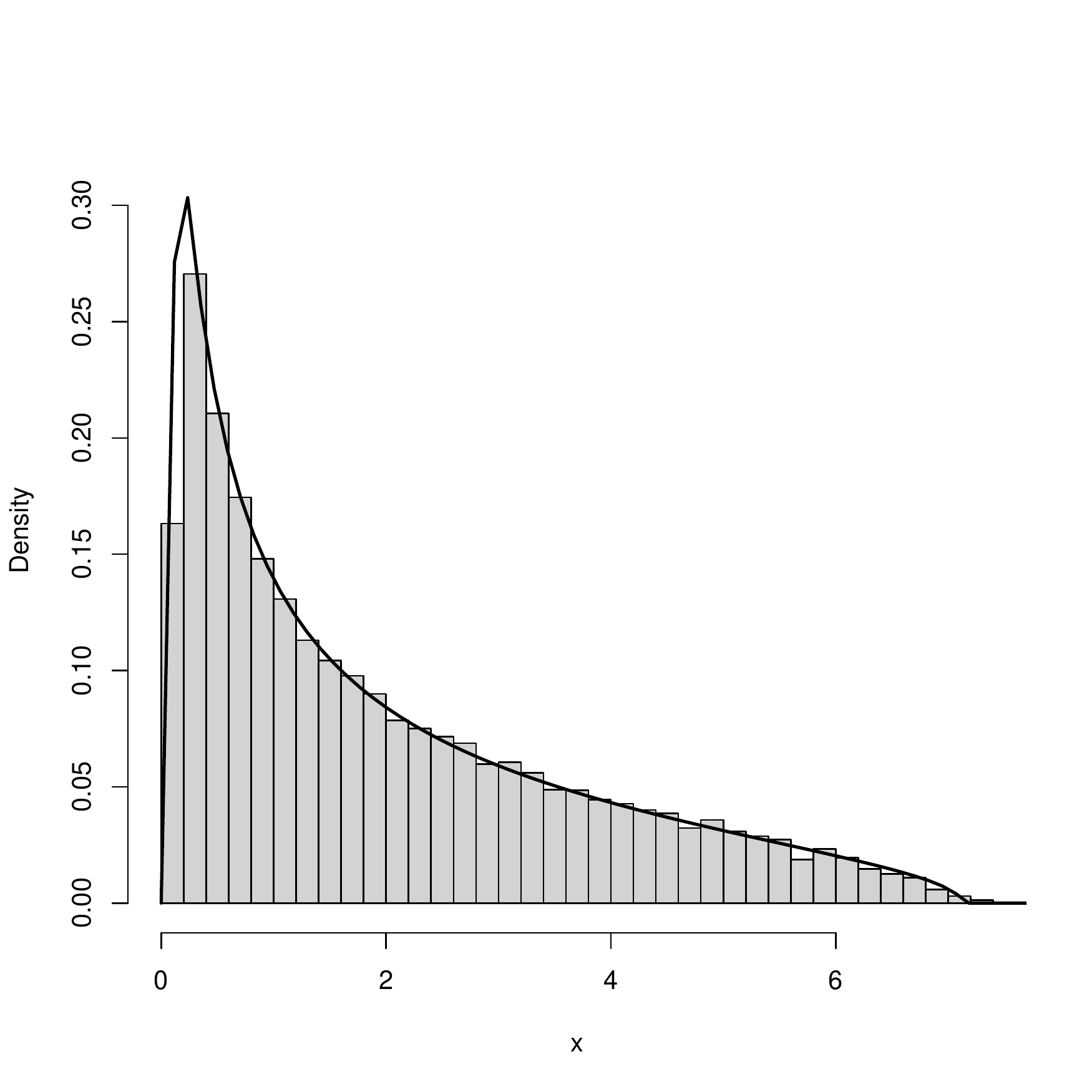,width=1.5 in,angle=0} &
			\psfig{figure=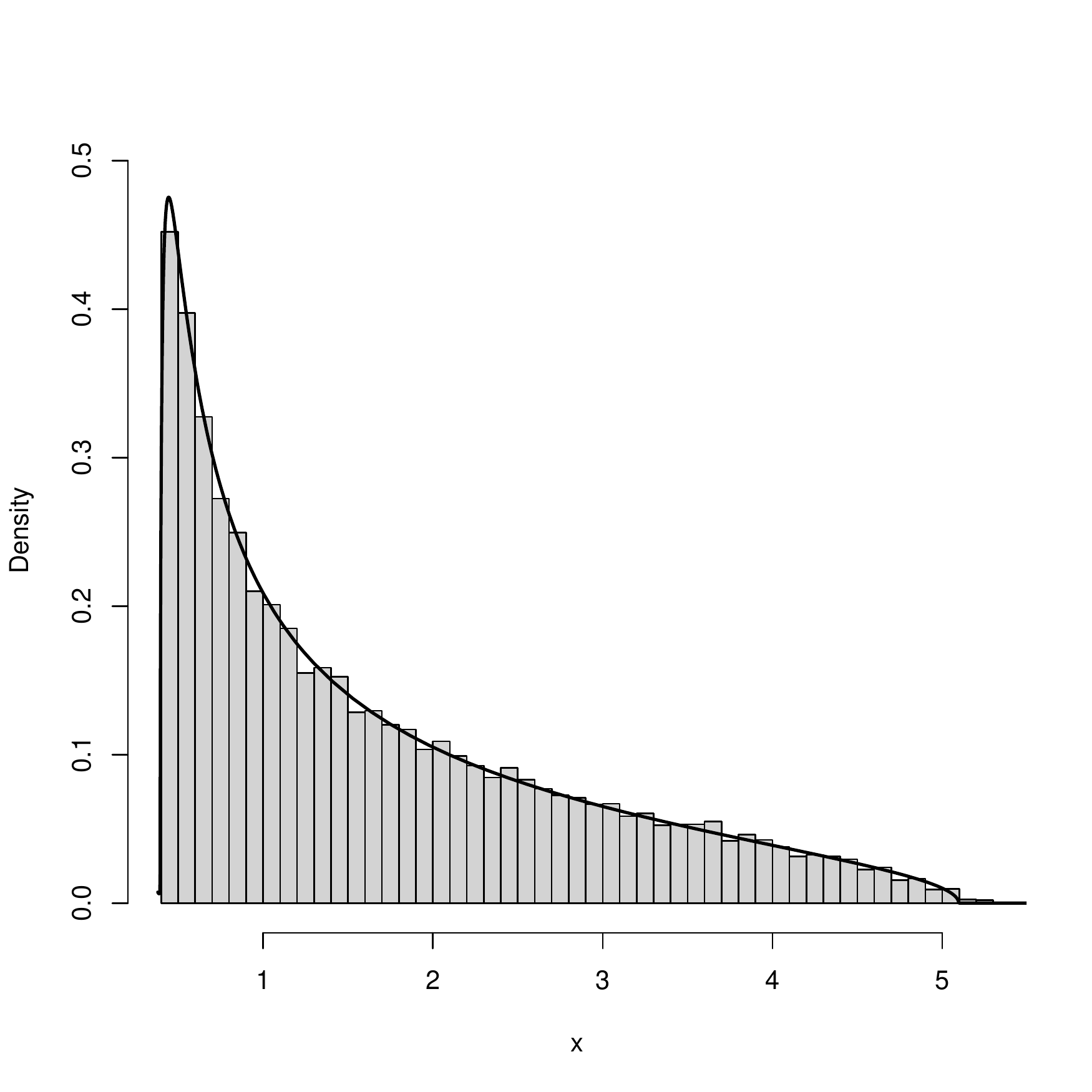,width=1.5 in,angle=0} \\
			\multicolumn{4}{c}{$(n,p)=(100,200)$}\\	
			Covariance& Pearson& Spearman& Kendall\\
		\end{tabular}
	}
	\caption{The limiting spectral distributions of sample covariance matrix, Pearson's correlation matrix, Spearman's correlation matrix and Kendall's correlation matrix where  $\X_1, \ldots, \X_n,~i.i.d. \sim N\left(0,\bSig(0.5)\right)$.}
	\label{fig2}
\end{figure}

For general $\bSig$, \citet{marvcenko1967distribution} derived the LSD of the sample covariance matrix whose Stieltjes transform $m$ is given by the Mar\u{c}enko-Pastur equation 
\begin{align*}
	m=\int \frac{1}{t(1-y-yzm)-z} dH(t),~H(t)=\lim_{p \to \infty} F^{\sbSig}(t).
\end{align*}
This distribution is called the generalized Mar\u{c}enko-Pastur law.  \citet{el2009concentration} proved that the Mar\u{c}enko-Pastur equation also holds for the Pearson's correlation matrix.  For the Spearman's correlation matrix, our Theorem \ref{thm1} shows that  
\begin{align*}
	\brho_{n} \eqlsd \tilde{\brho}_n \eqlsd \W_n= \frac{3}{n} \sum_{i=1}^n \A_i \A_i\trans,   
\end{align*}
where $\eqlsd$ denotes the two matrices share the same LSD. Our Theorem \ref{thm2} further proves that the LSD is also the generalized  Mar\u{c}enko-Pastur law whose Stieltjes transform $m$ is
\begin{align*}
	m=\int \frac{1}{t(1-y-yzm)-z} dH(t),~H(t)=\lim_{p \to \infty} F^{3\sbSig_2}(t).
\end{align*}
For the Kendall's correlation matrix, \citet{li2021eigenvalues} established that
\begin{align*}
	\btau_n \eqlsd  \frac{2}{n} \sum_{i=1}^n \A_i \A_i\trans+\bSig_{3},
\end{align*} 
and they presented the LSD in two equations where the LSD is not the generalized  Mar\u{c}enko-Pastur law anymore. Here
\begin{align*}
	\bSig_1=\frac{2}{\pi}\arcsin(\bSig), ~\bSig_2=\frac{2}{\pi}\arcsin(\bSig/2),~\bSig_{3}=\bSig_1-2 \bSig_{2}.
\end{align*}
In particular, we consider a specific tridiagonal covariance matrix
\begin{align*}
	\bSig(\rho)=\begin{pmatrix}
		1 & \rho &  &&\\
		\rho & 1 & \rho&&\\
		&\ddots&\ddots&\ddots&\\
		&&\rho&1 &\rho\\
		&&&\rho&1
	\end{pmatrix},
\end{align*}
where $\rho \neq 0$ and $|\rho| \leq 1/2$. By Szeg\"{o}  Theorem, we have
\begin{align*}
	\lim_{p \to \infty} F^{\sbSig}(t)=1-\frac{1}{\pi} \arccos \frac{t-1}{2 \rho},~|t-1|\leq 2|\rho|,
\end{align*}
and 
\begin{align*}
	\lim_{p \to \infty} F^{3\sbSig_2}(t)=1-\frac{1}{\pi} \arccos \frac{t-1}{2 \rho_1},~|t-1|\leq 2|\rho_1|,~\rho_1=\frac{6}{\pi} \arcsin\frac{\rho}{2}.
\end{align*}
Therefore, the Mar\u{c}enko-Pastur equation of the sample covariance matrix and the Pearson's correlation matrix is
\begin{align*}
	m=\int_{1-2|\rho|}^{1+2|\rho|} \frac{1}{t(1-y-yzm)-z} d\left(1-\frac{1}{\pi} \arccos \frac{t-1}{2 \rho}\right)=\frac{1}{\sqrt{(1-y-yzm-z)^2-4\rho^2(1-y-yzm)^2}},
\end{align*}
and similarly, the Mar\u{c}enko-Pastur equation of the Spearman's correlation matrix is
\begin{align*}
	m=\frac{1}{\sqrt{(1-y-yzm-z)^2-4\rho_1^2(1-y-yzm)^2}}.
\end{align*}
Solving the equation and using the inversion formula of the Stiejtjes transformation, we can get the LSDs. The LSD of Kendall's correlation matrix can be found in the Proposition 4.2 of \citet{li2021eigenvalues}. Fig \ref{fig2} shows the theoretical LSDs and the empirical distributions based on 100 replications for $\bSig=\bSig(0.5)$.

In summary, under mild conditions, the Pearson's correlation matrix shares the same properties of the sample covariance matrix. For rank-based correlation matrix, the Spearman's rho has the generalized  Mar\u{c}enko-Pastur law with the population covariance $6/\pi \arcsin(\bSig/2)$ and the Kendall's tau is equivalent to a sample covariance matrix with the population covariance matrix $\bSig_{2}$ plus a deterministic covariance matrix $\bSig_{3}$.

\section*{Acknowledgments}
We thank the Editor, an Associate Editor, and anonymous reviewers for their insightful comments. Wang's research is supported by National Natural Science Foundation of  China (11971017, 12031005) and NSF of Shanghai (21ZR1432900).

\section{Appendix}
We first collect several important lemmas. 
\begin{lemma}[Grothendieck's identity] \label{clalem1}
	Consider a bi-variate normal distribution:
	\begin{align*}
		\begin{pmatrix}
			z_1\\
			z_2
		\end{pmatrix} \sim N \Big\{\begin{pmatrix}
			0\\
			0
		\end{pmatrix}, \begin{pmatrix}
			1& \rho \\
			\rho &1
		\end{pmatrix} \Large\Big\},
	\end{align*}
	where $\rho \in [-1,1]$. Then, $\E\{\sign(z_1)\sign(z_2)\}=4 \E \{I(z_1,z_2>0)\}-1= 2 \arcsin{\rho}/\pi$.
\end{lemma}

\begin{lemma}[\citet{esscher1924method}] \label{clalem2}
	Consider a multivariate normal distribution 
	\begin{align*}
		\begin{pmatrix}
			z_1\\
			z_2\\
			z_3\\
			z_4
		\end{pmatrix} \sim N \left(\begin{pmatrix}
			0\\
			0\\
			0\\
			0
		\end{pmatrix},   \begin{pmatrix}
			1 & 1/2&\rho&\rho/2\\
			1/2&1&\rho/2&\rho\\
			\rho&\rho/2&1&1/2\\
			\rho/2&\rho&1/2&1
		\end{pmatrix} \right)
	\end{align*}
	where $\rho \in (-1,1)$. Then,
	\begin{align} \label{case0}
		\E \prod_{j=1}^4 \sign(z_j)=(\frac{2}{\pi} \arcsin{\rho})^2-(\frac{2}{\pi} \arcsin{\rho/2})^2+\frac{1}{9}.
	\end{align}
\end{lemma}

\begin{lemma} \label{lem-asin}
	For any $x \in [0,1]$, $2 \arcsin(x/2) \leq \arcsin(x) \leq 3 \arcsin(x/2)$ 
	and		$2 x/\pi \leq 2 \arcsin(x)/\pi \leq x$.
\end{lemma}

\begin{lemma} \label{aprop1}
	Assuming $\X_1, \ldots, \X_n,~ i.i.d.\sim N(\mathbf{0},\bSig)$, where $\Sigma_{ii}=1$,  we have $\E(\A_{ij})=\E(\A_i)=\mathbf{0}$ and
	\begin{align*}
		\cov(\A_{ij})=\bSig_1=\frac{2}{\pi}\arcsin(\bSig),~\cov(\A_i)=\bSig_2=\frac{2}{\pi}\arcsin(\bSig/2).
	\end{align*}
\end{lemma}

\begin{lemma}\label{lem2-var}
	Assume $\X_1, \ldots, \X_n,~ i.i.d.\sim N(\mathbf{0},\bSig)$, where $\Sigma_{ii}=1$, we have
	\begin{align*}
		\var( \A_{12} \trans \A_{13} ) \geq \var(\A_{12} \trans \A_{1}) \geq \var(\A_{1} \trans \A_{1})
	\end{align*}
	and 
	\begin{align}
		\var( \A_{12} \trans \A_{13} )=\tr(\bSig^2_1)-\tr(\bSig^2_2).
	\end{align}
\end{lemma}

\begin{myproof}{Lemma \ref{lem2-var}}
	Firstly, it is easy to see 	$\E (\A_{12} \A_{13} \trans)=\E (\A_{12} \A_{1} \trans)=\E (\A_{1} \A_{1} \trans)=\bSig_2$,
	which yields $\E (\A_{12} \trans  \A_{13}) =\E (\A_{12} \trans \A_{1}) =\E (\A_{1} \trans \A_{1}) =\tr(\bSig_2)=p/3$. 	To derive the inequality of variances, we use the Hoeffding decomposition
	\begin{align*}
		\var(\A_{12} \trans \A_{13} )=\E\{\A_{12} \trans \A_{13}-\tr(\bSig_2)\}^2
		=&\E\{\A_{12} \trans \A_{13}-\A_{12} \trans \A_{1}+\A_{12} \trans \A_{1}-\tr(\bSig_2)\}^2\\
		=&~\E(\A_{12} \trans \A_{13}-\A_{12} \trans \A_{1})^2+\var(\A_{12} \trans \A_{1})+2 \E(\A_{12} \trans \A_{13}-\A_{12} \trans \A_{1}){ (\A_{12} \trans \A_{1}-\tr(\bSig_2))}\\
		=&~\E(\A_{12} \trans \A_{13}-\A_{12} \trans \A_{1})^2+\var(\A_{12} \trans \A_{1}) \geq \var(\A_{12} \trans \A_{1}), 
	\end{align*}
	and similarly 
	\begin{align*}
		\var(\A_{12} \trans \A_{1} )=&~\E(\A_{12} \trans \A_{1}-\A_{1} \trans \A_{1})^2+\var(\A_{1} \trans \A_{1})+2 \E(\A_{12} \trans \A_{1}-\A_{1} \trans \A_{1})(\A_{1} \trans \A_{1}-\tr(\bSig_2))\\
		=&~\E(\A_{12} \trans \A_{1}-\A_{1} \trans \A_{1})^2+\var(\A_{1} \trans \A_{1}) \geq \var(\A_{1} \trans \A_{1}). 
	\end{align*}
The calculation of $\var( \A_{12} \trans \A_{13} )$ can be found in the proof of Lemma A 3.2 of \citet{li2021eigenvalues}. The proof of Lemma~\ref{lem2-var} is completed.  
\end{myproof}

\begin{myproof}{Proposition \ref{prob2}}
Here we prove Proposition \ref{prob2} and Proposition \ref{prob1} is a special case of Proposition \ref{prob2}  by setting $\bSig=\bI$.

Noting $\R \trans=(\R_1, \ldots,\R_n)$, we have
\begin{align*}
	\E \{\frac{1}{n} \R \trans \R\}=\frac{1}{n} \sum_{i=1}^n \R_i \R_i\trans=\E \{ \R_1\R_1\trans\},~\E \{\frac{1}{p} \R  \R \trans\}=\frac{1}{p}\E \left( \R_i\trans \R_j \right)_{n \times n}.
\end{align*}

By the formula \eqref{Ri}, 
\begin{align*} 
	\R_i=\sqrt{\frac{3}{n^2-1}}\sum_{k\neq i}\A_{ik},
\end{align*}
and we have
\begin{align*}
	\E(\R_1\R_1\trans)=\frac{3}{n^2-1} \sum_{l,m=2}^n \E \{\A_{1l}\A_{1m}\trans \}
	=&\frac{3}{n^2-1}\{(n-1) \cov(\A_{12})+(n-1)(n-2) \cov(\A_1) \} \\
	=&\frac{3}{n+1}\left(\bSig_1+(n-2) \bSig_{2} \right).
\end{align*}

Next, we calculate $\E(\R_i\trans \R_j)$. When $j=i$, 
\begin{align*}
	\E(\R_i\trans \R_i)=\tr\left\{ 	\E(\R_1\R_1\trans) \right\}=\frac{3}{n+1} \left(\tr(\bSig_1)+(n-2) \tr(\bSig_{2}) \right)=\frac{3}{n+1}\left( p+\frac{n-2}{3}p\right)=p.
\end{align*}
For $i\neq j$,
\begin{align*}
	\E(\R_i\trans \R_j)&=\frac{3}{n^2-1}\E\left(\sum_{l\neq i}\A_{il}\right)\trans\left(\sum_{m\neq j}\A_{jm}\right)=\frac{3}{n^2-1}\E\sum_{l \neq i,m \neq j} \A_{il}\trans \A_{jm}\\
	&=\frac{3}{n^2-1}\left( \E\sum_{m\neq j}\A_{ij}\trans \A_{jm}+ \E\sum_{m\neq i,j}\A_{im}\trans \A_{jm}+\E \sum_{l\neq i, l \neq j}\A_{il}\trans \A_{ji}\right)\\
	&=\frac{3}{n^2-1}\left( -(n-2)\tr(\bSig_2)-\tr(\bSig_1)+(n-2)\tr(\bSig_2)-(n-2)\tr(\bSig_2) \right)=-\frac{3}{n^2-1}\left(p+\frac{n-2}{3}p\right)=-\frac{p}{n-1}.
\end{align*}
Thus
\begin{align*}
	\E \{\frac{1}{p} \R  \R \trans\}=\frac{n}{n-1}(\bI_n-\frac{1}{n}\one_n\one_n\trans).
\end{align*}
The proof  is completed.  
\end{myproof}

\begin{myproof}{Theorem \ref{thm1}} Recall that 
\begin{align*}
	\tilde{\brho}_{n}-\W_n=\frac{3}{n(n-1)(n-2)} \sum_{i,j,k}^* \left(\A_{ij}\A_{ik}\trans-\A_i\A_i\trans\right).
\end{align*}
For $i,j,k$, defining the kernel function $\phi(i,j,k)=\A_{ij}\A_{ik}\trans -\A_i\A_i\trans$ and the symmetric kernel function
\begin{align*}
	\psi(i,j,k)=\frac{1}{6} \sum_{(i',j',k')=\pi(i,j,k)}	\phi(i',j',k'),
\end{align*}
we have
\begin{align*}
	\tilde{\brho}_{n}-\W_n=\frac{18}{n(n-1)(n-2)}  \sum_{i<j<k} \psi(i,j,k).
\end{align*}
By the symmetric properties of U-statistics, we can get
\begin{align} \label{thm2p0}
	\frac{1}{9}\E \{ \left( 	\tilde{\brho}_{n}-\W_n \right)  \trans \left( 	\tilde{\brho}_{n}-\W_n \right)  \}
	=&\frac{6}{n(n-1)(n-2)}  \sum_{i<j<k} \E\{\psi(1,2,3) \trans \psi(i,j,k) \} \nonumber \\
	=&\frac{6}{n(n-1)(n-2)}  \E\{\psi(1,2,3) \trans \psi(1,2,3) \} +\frac{18(n-3)}{n(n-1)(n-2)}  \E\{\psi(1,2,3) \trans \psi(1,2,4) \} \nonumber \\
	&+\frac{18(n-3)(n-4)}{2n(n-1)(n-2)}  \E\{\psi(1,2,3) \trans \psi(1,4,5) \}.
\end{align}
Next, we bound these terms. For the first term, we have
	\begin{align*}
	\E \left\|\phi(1,2,3)\right\|_2^2= \E \tr \left( \A_{12}\A_{13}\trans-\A_1\A_1\trans  \right) \left( \A_{13}\A_{12}\trans-\A_1\A_1\trans  \right)=\E \{(\A_{12} \trans \A_{12})  (\A_{13} \trans \A_{13})  \}-\E (\A_{1} \trans \A_{1})^2 \leq p^2
	\end{align*}
	and then
	\begin{align}
	\E \|\psi(1,2,3) \|_2^2= \E \left\|\frac{1}{6} \sum_{(i,j,k)=\pi(1,2,3)}	\phi(i,j,k)\right\|_2^2 
		\leq \frac{1}{6} \sum_{(i,j,k)=\pi(1,2,3)}	\E \left\| \phi(i,j,k)\right\|_2^2= \E \left\|\phi(1,2,3)\right\|_2^2\leq p^2.  		  \label{thm2p1}
	\end{align}
For the second term, 
	\begin{align} \label{thm2p2}
		\E\{\psi(1,2,3) \trans \psi(1,2,4) \} \leq&  \left\{   \E\| \psi(1,2,3)\|_2^2 \right\}^{1/2}  \left\{   \E\| \psi(1,2,4)\|_2^2 \right\}^{1/2} =\E \|\psi(1,2,3) \|_2^2 \leq   p^2.	  
	\end{align}
 For the third term, we have the conditional expectation
	\begin{gather*}
		\E\{\phi(1,2,3)\mid \X_1\}=	\E\{\phi(1,3,2)\mid \X_1\}=\bzero,~\E\{\phi(2,1,3)\mid \X_1\}=	\E\{\phi(3,1,2)\mid \X_1\}=\E\{\A_{21}\A_2\trans-\bSig_{2} \mid \X_1\},\\
		\E\{\phi(2,3,1)\mid \X_1\}=	\E\{\phi(3,2,1)\mid \X_1\}=\E\{\A_{2}\A_{21}\trans-\bSig_{2} \mid \X_1\},
	\end{gather*} 
	and then
	\begin{align*}
		& \E\{\psi(1,2,3)\mid \X_1\}=\frac{1}{3}\E\{(\A_{2}\A_{21}\trans-\bSig_{2})+(\A_{21}\A_{2}\trans-\bSig_{2}) \mid \X_1\}.
	\end{align*}
	Thus,
	\begin{align*}
	\E\{\psi(1,2,3) \trans \psi(1,4,5) \}=& \frac{1}{9}\E\tr \{(\A_{2}\A_{21}\trans-\bSig_{2})+(\A_{21}\A_{2}\trans-\bSig_{2})\} \{(\A_{3}\A_{31}\trans-\bSig_{2})+(\A_{31}\A_{3}\trans-\bSig_{2})\}\\
		=&\frac{2}{9}\{\E (\A_{21}\trans \A_3)(\A_{31} \trans \A_2)+\E(\A_{21}\trans \A_{31})(\A_2\trans \A_3)-4\tr(\bSig^2_{2})  \}\\
		\leq &\frac{1}{9}\{\E (\A_{21}\trans \A_3)^2+\E(\A_{31} \trans \A_2)^2+2 \E(\A_{21}\trans \A_{31}-\tr(\bSig_{2}))(\A_2\trans \A_3) \}-\frac{8}{9}  \tr(\bSig^2_{2}) \\
		\leq & \frac{1}{9}\{\E (\A_{21}\trans \A_3)^2+\E(\A_{31} \trans \A_2)^2\}+\frac{2}{9} \{\var(\A_{12}\trans \A_{13}) \}^{1/2} \{\var(\A_{2}\trans \A_3) \}^{1/2}-\frac{8}{9}  \tr(\bSig^2_{2})\\
		=& \frac{2}{9} \tr(\bSig_1 \bSig_{2})+\frac{2}{9} \{\tr(\bSig^2_1)-\tr(\bSig^2_2)\}^{1/2} \{ \tr(\bSig^2_{2})\}^{1/2}-\frac{8}{9}  \tr(\bSig^2_{2}).
	\end{align*}
	Together with Lemmas  \ref{lem-asin} and \ref{lem2-var}, we can get
	\begin{align} \label{thm2p3}
		\E\{\psi(1,2,3) \trans \psi(1,4,5) \} \leq \{ \frac{6}{9} +\frac{2\sqrt{8}}{9}-\frac{8}{9}\}  \tr(\bSig^2_{2})\leq \frac{1}{2} \tr(\bSig^2_{2}).
	\end{align}

Finally, combining \eqref{thm2p0}- \eqref{thm2p3}, we conclude that 
\begin{align*}
	\frac{1}{9}\E \{ \left( 	\tilde{\brho}_{n}-\W_n \right)  \trans \left( 	\tilde{\brho}_{n}-\W_n \right)  \} 
	\leq  \frac{6+18(n-3)}{n(n-1)(n-2)} p^2+\frac{9(n-3)(n-4)}{n(n-1)(n-2)} \{ \tr(\bSig^2_{2})\}.
\end{align*}
By Corollary A.41 of \citet{bai2010spectral}, 
\begin{align*}
	L^3\left(F^{	\tilde{\brho}_{n}}, F^{\W_n} \right) \leq \frac{1}{p} \|\tilde{\brho}_{n}-\W_n\|_2^2 
\end{align*}
which yields $\E 	L^3\left(F^{	\tilde{\brho}_n}, F^{\W_n} \right)  \to 0$ and then 
\begin{align*}
	L\left(F^{\tilde{\brho}_n}, F^{\W_n} \right)  \to 0,~\text{in probability}.
\end{align*}
The proof  is completed.  
\end{myproof}

\begin{myproof}{Proposition \ref{prop4}} Since 
\begin{align*}
	\brho_n=\frac{3}{n+1} \btau_n-\frac{3}{n+1} \tilde{\brho}_n+\tilde{\brho}_n,
\end{align*}
by Corollary A.41 of \citet{bai2010spectral}, 
\begin{align*}
	L^3\left(F^{\brho}_{n}, F^{\tilde{\brho}_{n}} \right) \leq \frac{1}{p} \|\brho_n-\tilde{\brho}_{n}\|_2^2 \leq \frac{18}{n^2 p} \|\btau_n\|_2^2+ \frac{18}{n^2 p} \|\tilde{\brho}_n\|_2^2.
\end{align*}
Noting $\btau_n$ is a correlation matrix, we have $\| \btau_n\|_\infty \leq 1$ and then 
$\| \btau_n\|_2^2 \leq p^2$.

Similarly, we can show $\|\tilde{\brho}_{n}\|_2^2 \leq 9 p^2$. Thus  $L\left(F^{\brho_n}, F^{\tilde{\brho}_{n}} \right) \to 0$. Together with Theorem \ref{thm1}, the claim follows.   
\end{myproof}

\begin{myproof}{Theorem \ref{thm1-1}}
Recall
\begin{align*}
	\R=\sqrt{\frac{12}{n^2-1}} \left(r_{ij}-\frac{n+1}{2}\right)_{n \times p}=  \sqrt{\frac{12n^2}{n^2-1}}   \left( \hat{F}_j(x_{ij})-\frac{n+1}{2n} \right)_{n \times p}
\end{align*}
and $\brho_{n}=\frac{1}{n} \R \trans \R$. Writing
\begin{align*}
	\tilde{\R}=\sqrt{\frac{12n^2}{n^2-1}}   \left( F_j(x_{ij})-\frac{n+1}{2n} \right)_{n \times p},
\end{align*}
by Corollary A.42 of \citet{bai2010spectral}, we have
\begin{align*}
	L^4\left(F^{\brho_n}, F^{\frac{1}{n} \tilde{\R} \trans \tilde{\R}} \right) \leq  \frac{2}{p^2} \left\{ \tr \left(\brho_n+ \frac{1}{n} \tilde{\R} \trans \tilde{\R}\right) \right\} \left\{ \frac{1}{n} \|\R-\tilde{\R}\|_2^2 \right\}
	\leq  \frac{20}{np} \|\R-\tilde{\R}\|_2^2, 
\end{align*}
where we use the facts $\tr(\brho_{n})=p$ and 
\begin{align*}
	\tr \left(\frac{1}{n} \tilde{\R} \trans \tilde{\R}\right) \leq p \| \tilde{\R}\|^2_\infty\leq  \frac{12n^2 p}{n^2-1}  \left(\frac{n+1}{2n} \right)^2 \leq 9p, ~(n \geq 2).
\end{align*}
Noting
\begin{align*}
	\frac{1}{np} \|\R-\tilde{\R}\|_2^2\leq \frac{12n^2}{n^2-1} \max_{i,j} |\hat{F}_j(x_{ij})-F_j(x_{ij}) |^2
	\leq  \frac{12n^2}{n^2-1} \max_{j} \left\{ \sup_{x \in \mR} |\hat{F}_j(x)-F_j(x)|^2 \right\}   
	= \frac{12n^2}{n^2-1}  \left\{\max_{j} \|\hat{F}_j(x)-F_j(x)\|_\infty \right\}^2,
\end{align*}
we need to control $\|\hat{F}_j(x)-F_j(x)\|_\infty,~j \in \{1,\ldots, p\}$ simultaneously. Thus, we shall use the Dvoretzky–Kiefer–Wolfowitz (DKW) inequality which provides a refined result of the classical Glivenko–Cantelli Theorem. By DKW inequality, $\forall \epsilon>0$,
\begin{align*}
	\Pr \left( \|\hat{F}_j(x)-F_j(x)\|_\infty \geq \epsilon  \right) \leq 2 e^{-2 n\epsilon^2}
\end{align*} 
which yields 
\begin{align*}
	\Pr \left(\max_{j} \|\hat{F}_j(x)-F_j(x)\|_\infty \geq \epsilon  \right) \leq \sum_{j=1}^p \Pr \left( \|\hat{F}_j(x)-F_j(x)\|_\infty \geq \epsilon  \right) \leq 2 p e^{-2 n\epsilon^2}.
\end{align*}
Since $p/n \to y$, we have $2 p e^{- n\epsilon^2} \leq 1 $ for large enough $n$ and then 
\begin{align*}
	\sum_{n=1}^\infty 2 p e^{-2 n\epsilon^2}=\sum_{n=1}^\infty \{2 p e^{- n\epsilon^2}\}\cdot e^{- n\epsilon^2} < \infty.
\end{align*}
By the Borel–Cantelli Lemma, we obtain 
\begin{align*}
	\max_{j} \|\hat{F}_j(x)-F_j(x)\|_\infty \to 0, ~\text{almost surely},
\end{align*} 
and then
\begin{align} \label{as1}
	L\left(F^{\brho_n}, F^{\frac{1}{n} \tilde{\R} \trans \tilde{\R}} \right) \to 0,~\text{almost surely}.
\end{align}

Next, we bound the difference between $ \tilde{\R} \trans \tilde{\R}/n$ and $\W_n=3\sum_{i=1}^n \A_i \A_i\trans/n$. Noting
\begin{align*}
	\tilde{\R}=\sqrt{\frac{12n^2}{n^2-1}}   \left( F_j(x_{ij})-\frac{n+1}{2n} \right)_{n \times p}
\end{align*}
and 
\begin{align*}
	\sqrt{3}	(\A_1,\ldots,\A_n)\trans= \sqrt{12} \left(F_j(x_{ij})-\frac{1}{2}\right)_{n \times p},
\end{align*}
we can get
\begin{align*}
	\left\|	\tilde{\R}- \sqrt{3}	(\A_1,\ldots,\A_n)\trans\right\|_\infty\leq 2 \sqrt{3} \left|\sqrt{\frac{n^2}{n^2-1}} -1\right|+\sqrt{3}  \left|\frac{n+1}{\sqrt{n^2-1} } -1\right| \to 0.
\end{align*}
Using Corollary A.42 of \citet{bai2010spectral} again, we have $L^4\left(F^{\W_n}, F^{\frac{1}{n} \tilde{\R} \trans \tilde{\R}} \right) \to 0$. Combined with \eqref{as1}, we conclude that 
\begin{align*} 
	L\left(F^{\brho_n}, F^{\W_n} \right) \to 0, ~\text{almost surely}.
\end{align*}
The proof  is completed.  
\end{myproof}

\begin{myproof}{Theorem \ref{thm2}}
Since $\A_1,\ldots, \A_n$ are i.i.d., we can apply Theorem 1.1 of \citet{bai2008large}. It is sufficient to check assumption 1 of their theorem, i.e.,  for any non-random $p\times p$ matrix $\mathbf{B}$, 
\begin{align*}
	\frac{1}{p^2}	\var\left(\A_1\trans \mathbf{B} \A_1 \right) \to 0.
\end{align*}
Noting $\A_i=2 \Phi(\X_i)-1$ where $\Phi(\cdot)$ is differentiable,  the Gaussian Poincar\'{e} inequality together with $\|\bSig\| \leq C$ can yield the above conclusion. A detailed proof can be found in Lemma 3.3 of \cite{li2021eigenvalues}.  The proof  is completed.  
\end{myproof}

\bibliography{ref}
\end{document}